\documentclass[11pt,a4paper]{amsart}

\usepackage{times,epsfig,epic,eepic}
\usepackage{amsfonts,amscd,amsmath,amssymb,amsthm}

\copyrightinfo{2004}{by the authors}

\numberwithin{equation}{section}

\renewcommand{\epsilon}{\varepsilon}

\DeclareSymbolFont{SY}{U}{psy}{m}{n}
\DeclareMathSymbol{\emptyset}{\mathord}{SY}{'306}

\renewcommand{\i}{\mathrm{i}}
\newcommand{\e}{\mathrm{e}}

\newcommand{\sk}{\mathsf{k}}

\newcommand{\bw}{\mathbf{w}}

\newcommand{\fL}{\mathfrak{L}}

\newtheorem{theorem}{Theorem}[section]{\bf}{\it}
\newtheorem{proposition}[theorem]{Proposition}{\bf}{\it}
\newtheorem{corollary}[theorem]{Corollary}{\bf}{\it}
\newtheorem{example}[theorem]{Example}{\it}{\rm}
\newtheorem{lemma}[theorem]{Lemma}{\bf}{\it}
\newtheorem{remark}[theorem]{Remark}{\it}{\rm}
\newtheorem{definition}[theorem]{Definition}{\bf}{\it}
{\bf}{\it}

\DeclareMathOperator{\Ran}{Ran} \DeclareMathOperator{\Ker}{Ker}

\newcommand{\R}{\mathbb{R}}
\newcommand{\C}{\mathbb{C}}
\newcommand{\D}{\mathbb{D}}
\newcommand{\T}{\mathbb{T}}
\newcommand{\Z}{\mathbb{Z}}
\newcommand{\N}{\mathbb{N}}

\newcommand{\1}{\mathbb{I}}

\newcommand{\fH}{\mathfrak{H}}

\newcommand{\cA}{{\mathcal A}}
\newcommand{\cB}{{\mathcal B}}

\newcommand{\cD}{{\mathcal D}}
\newcommand{\cE}{{\mathcal E}}

\newcommand{\cG}{{\mathcal G}}
\newcommand{\cH}{{\mathcal H}}
\newcommand{\cI}{{\mathcal I}}

\newcommand{\cK}{{\mathcal K}}
\newcommand{\cL}{{\mathcal L}}
\newcommand{\cM}{{\mathcal M}}
\newcommand{\cN}{{\mathcal N}}

\newcommand{\cP}{{\mathcal P}}

\newcommand{\cS}{{\mathcal S}}

\newcommand{\cV}{{\mathcal V}}

\newcommand{\cW}{{\mathcal W}}

\renewcommand{\Im}{{\ensuremath{\mathrm{Im}\,}}}
\renewcommand{\Re}{{\ensuremath{\mathrm{Re}\,}}}

\begin{document}

\title[Random Walks on Graphs]{Generating Functions of Random Walks on Graphs}

\author[V.~Kostrykin]{Vadim Kostrykin}
\address{Vadim Kostrykin\\ Fraunhofer-Institut f\"{u}r Lasertechnik, Steinbachstra{\ss}e
15, D-52074\\ Aachen, Germany}
\email{kostrykin@ilt.fraunhofer.de, kostrykin@t-online.de}

\author[R.~Schrader]{Robert Schrader}
\address{Robert Schrader\\ Institut f\"{u}r Theoretische Physik\\
Freie Universit\"{a}t Berlin, Arnimallee 14\\ D-14195 Berlin, Germany}
\email{schrader@physik.fu-berlin.de}
\thanks{The work of R.S.\ was supported in part by DFG SFB
288 ``Differentialgeometrie und Quantenphysik''}

\keywords{Random walks, differential operators on graphs}

\subjclass[2000]{05A15, 34B45, 60G50}

\dedicatory{Dedicated to Ludwig Faddeev on the occasion of his 70th birthday\\
and to Philippe Blanchard and Konrad Osterwalder on the occasion
of their 60th birthday}

\begin{abstract}
The article provides an explicit algebraic expression for the generating
function of walks on graphs. Its proof is based on the scattering theory
for the differential Laplace operator on non-compact graphs.
\end{abstract}

\maketitle

\section{Introduction}\label{sec:1}

The concept of a generating function is known to be a very important tool in
combinatorics, probability, and number theory. Associated methods reduce
the solution of combinatorial or probabilistic problems to the study of
particular properties of the generating function which can be performed by
methods of function theory and analysis. For an introduction to this method
the reader may consult the books \cite{Egorychev}, \cite{Graham},
\cite{Wilf}. A number of solved and still unsolved combinatorial problems,
where the generating function plays a central role can be found in the
article \cite{Shapiro}. In the probabilistic context we mention the
solution of the problem whether a simple random walk on $\Z^d$ is recurrent
or transitive by an analysis of the generating function (see, e.g.,
\cite{Grinstead}). Some further examples will be discussed below in
Sections \ref{sec:2} and \ref{sec:random}.

The present work is devoted to the determination of the generating function
for walks on graphs (both in combinatorial and probabilistic contexts).
Walks on graphs are considered, in particular, in \cite{Aldous:Fill},
\cite{Doyle}, \cite{Sunada}, \cite[Section 4.7]{Stanley}, \cite{Woess}. In
a custom setting random walks on graphs are defined as Markov chains on the
vertices of the graph. The transition probability from one vertex to
another is assumed to be non-zero if and only if these vertices are
adjacent. For a survey of the theory of random walks on graphs, see
\cite{Lovasz}.

We consider a slightly different but closely related model of random walks
on graphs, where the states are chosen to be the edges of the graph.
Transitions between different states are determined by stochastic (that is,
Markov) matrices $M(v)$ prescribed at every vertex of the graph. The graphs
are assumed to be non-compact, that is, besides a finite number of edges
(or ``internal lines'' $i\in\cI$) to have a non-empty
set $\cE$ of ``external lines'' which serve as entries or exits for random
walks. More precisely the model will be described in the following section.
A relation between this model of random walks and random walks on vertices
is explained in the Appendix below.

Consider an arbitrary positive weight on the graph (that is, a map
assigning to any edge $i$ of the graph a positive number $a_i$). We will
call $\underline{a}=\{a_i\}_{i\in\cI}\in(\R_+)^{|\cI|}$ a \emph{penalty
vector}. With $\beta$ being a complex parameter we define a generating
function $T_{e,e^\prime}(\beta)$ of walks from an external line
$e^\prime\in\cE$ to an external line $e\in\cE$ as
\begin{equation*}
\begin{split}
T_{e,e^\prime}(\beta) & = \sum [M(v_N)]_{e,i_N} \e^{-\beta a_{i_N}}
[M(v_{N-1})]_{i_N,i_{N-1}}\ldots \\ & \qquad\qquad \ldots
[M(v_{1})]_{i_2,i_1} \e^{-\beta a_{i_1}} [M(v_{0})]_{i_1,e^\prime},
\end{split}
\end{equation*}
where the sum is taken over all walks
$\{e^\prime,i_1,i_2,\ldots,i_{N-1},i_N,e\}$ from $e^\prime$ to $e$, the set
$v_0,v_1,\ldots v_N$ is the ordered list of vertices (with possible
repetitions) visited by the walk, $i_k$ are the corresponding internal
lines traversed during the walk (again with possible repetitions). In the
context of the generating functions the weight $\exp\{-\beta a_i\}$ can be
viewed as a penalty factor for traversing the edge $i$ during a walk.

More generally, one can also consider penalty vectors depending on the
direction in which a given edge is traversed by the walk. The corresponding
generating function will be discussed in Section \ref{sec:3} below (see
Theorem \ref{genrep}).

The main result of the present work (see Theorem \ref{thm:3.2} below)
provides an explicit algebraic expression for the generating function of
walks on graphs. Its proof is based on the scattering theory for the
\emph{differential} Laplace operator on non-compact graphs and the
corresponding methods developed by the authors in \cite{KS1}, \cite{KS2},
\cite{KS3}, \cite{KS4}, \cite{KS5}. In the context of differential
operators the weights $a_i$ will be interpreted as the metric lengths of
the edges $i$.

The generating function is determined by analytic continuation of the
scattering matrix to complex values of the spectral parameter. This result
is very reminiscent of a similar result in relativistic quantum field theory
in the context of vacuum expectations of products of quantum fields. The
analytic continuation of the Wightman distributions \cite{Wightman} to the
Euclidean points (the so called Wick rotation \cite{Wick}) results in the
Schwinger functions \cite{Schwinger}. Conversely, by a result \cite{OS1},
\cite{OS2} of K.~Osterwalder and one of the authors (R.~S.) the Schwinger
functions give rise to Wightman distributions. In the bosonic case Symanzik
and Nelson have shown that the Schwinger functions describe a stochastic
theory (see \cite{GJ}, \cite{Nelson}, \cite{Symanzik}, and references
quoted there).

We expect that the model of walks on graphs considered in the present
article may be of interest in the context of optimization of traffic flows
and in telecommunication networks, where the transition matrices $M(v)$
determine a proportion of the traffic or signals to be transmitted in a
given direction.

The article is organized as follows. In Section \ref{sec:2} we will give
definitions of walks on graphs and of associated generating functions. Also
we present several examples relating the generating function to
combinatorics. In Section \ref{sec:operator} we will revisit the scattering
theory of differential Laplace operators on graphs. Section
\ref{sec:harmony} is devoted to the proof of the combinatorial Fourier
expansion formula \eqref{sseries}. Theorem \ref{thm:main:harmony} proves
the absolute convergence of the Fourier series and Theorem \ref{verloren}
expresses the Fourier coefficients as sums over the walks on the graph. In
Section \ref{sec:analyt} we will consider the analytic continuation of the
scattering matrix with respect to the square root of the energy (that is,
the spectral parameter). In Section \ref{sec:3} the generating function
will be expressed in terms of the scattering matrix for a Laplace operator
with boundary conditions determined by the transition matrices $M(v)$. In
Section \ref{sec:random} we will turn to random walks on graphs. By means
of the generating function we will calculate several mean values associated
to this probabilistic set-up.

There are several further models which can also be treated by the methods
of the present work. In particular, choosing the matrices $M(v)$ as
independent random variables one obtains a model of random walks in random
environment. Further, the graph itself can be chosen to be random (see,
e.g., \cite{Bollobas}). We note that random graphs have been used to model
the spread of epidemics like AIDS, see, e.g., the article \cite{BCK} and
further references quoted there.

\section{Walks on Graphs}\label{sec:2}

We consider a finite, connected and non-compact graph $\cG=(V,\cI,\cE,\partial)$, where
$V=V(\cG)$ is a finite set of \emph{vertices}, $\cI$ is a finite set of
\emph{internal lines}, $\cE$ is a finite set of \emph{external lines}. The
elements of the set $\cI\cup\cE$ are called \emph{edges}. The boundary
operator $\partial$ assigns to each internal line $i\in\cI$ an ordered pair
$(v_1,v_2)$ of vertices (possibly equal) and to each external line
$e\in\cE$ a single vertex $v$. The vertices $v_1:=\partial^-(i)$ and
$v_2:=\partial^+(i)$ are called the \emph{initial} and \emph{terminal}
vertex of the internal line $i$, respectively. This obviously induces
an orientation on each of the internal lines and this will become
relevant below.

The vertex $v=\partial(e)$ is the initial vertex of the external line $e$.
If $\partial(i)=(v,v)$, then $i$ is called a \emph{tadpole}. To simplify
the discussion, in what follows we will assume that the graph $\cG$
contains no tadpoles.

Two vertices $v$ and $v^\prime$ are called \emph{adjacent} if there is an
internal line $i\in\cI$ such that either $(v,v^\prime)=\partial(i)$ or
$(v^\prime,v)=\partial(i)$. A vertex $v$ and the (internal or external) line
$j\in\cI\cup\cE$ are \emph{incident} if $v\in\partial(j)$. The
\emph{degree}  $\deg(v)$ equals the number of (internal or external) lines
incident with the vertex $v$.

We do not require the map $\partial:\; \cI\rightarrow V\times V,\quad
\cE\rightarrow V$ to be one-to-one. In particular, any two vertices are
allowed to be adjacent to more than one internal line and two different
external lines may be incident with the same vertex.

Given an arbitrary vector $\underline{a}=\{a_i\}_{i\in\cI}\in\R^{|\cI|}$
with strictly positive components, we will endow the graph with the
following metric structure. Any internal line $i\in\cI$ will be associated
with an interval $[0,a_i]$ with $a_i>0$ such that the initial vertex of $i$
corresponds to $x=0$ and the terminal one - to $x=a_i$. Any external line
$e\in\cE$ will be associated with a half-line $[0,+\infty)$. The number
$a_i$ can be viewed as the length of the internal line $i$.

A nontrivial \emph{walk} $\bw$ on the graph $\cG$ from $e^\prime\in\cE$ to
$e\in\cE$ is a sequence
\begin{equation*}
\{e^\prime,i_1,\ldots,i_N,e\}
\end{equation*}
of edges such that
\begin{itemize}
\item[(i)]{$v_0 := \partial(e^\prime)\in\partial(i_1)$, $v_N :=
\partial(e)\in\partial(i_N)$, and for any $k\in\{1,\ldots,N-1\}$ there is a
vertex $v_k\in V$ such that $v_k\in\partial(i_k)$ and
$v_k\in\partial(i_{k+1})$};
\item[(ii)]{$v_k\neq v_{k+1}$ for all
$k\in\{0,\ldots,N-1\}$.}
\end{itemize}
The number $N$ is the \emph{combinatorial length}
$|\bw|_{\mathrm{comb}}\in\N$ and the number
\begin{equation*}
|\bw| = \sum_{k=1}^N a_{i_k} > 0
\end{equation*}
is the \emph{metric length} of the walk $\bw$.

\begin{example}
Let $\cG=(V,\cI,\cE,\partial)$ with $V=\{v_0,v_1\}$, $\cI=\{i\}$,
$\cE=\{e\}$, $\partial(e)=v_0$, and $\partial(i)=(v_0,v_1)$. Then the
sequence \{e,i,e\} is not a walk, whereas \{e,i,i,e\} is a walk from $e$ to
$e$.
\end{example}

\begin{proposition}\label{prop:2.1}
Given an arbitrary nontrivial walk $\bw=\{e^\prime,i_1,\ldots,i_N,e\}$
there is a unique sequence $\{v_k\}_{k=0}^N$ of vertices such that $v_0 =
\partial(e^\prime)\in\partial(i_1)$, $v_N =
\partial(e)\in\partial(i_N)$, $v_k\in\partial(i_k)$, and
$v_k\in\partial(i_{k+1})$.
\end{proposition}

\begin{proof}
Assume on the contrary that there are two different sequences
$\{v_k\}_{k=0}^N$ and $\{v^\prime_k\}_{k=0}^N$ satisfying the assumption of
the proposition. This implies that there is a number $K\in\{0,\ldots,N-2\}$
such that $v_k=v^\prime_k$ for all $k\in\{0,\ldots,K\}$ but $v_{K+1}\neq
v^\prime_{K+1}$. Obviously, the vertices $v_K$, $v_{K+1}$, and
$v^\prime_{K+1}$ are incident with the same edge. Thus, either
$v_K=v_{K+1}$ or $v_K=v^\prime_{K+1}$, which is a contradiction.
\end{proof}

We emphasize, that at any vertex of the sequence $\{v_k\}_{k=0}^N$
associated with a nontrivial walk $\bw$, the walk is either ``reflected''
or ``transmitted''.

A \emph{trivial} walk $\bw$ on the graph $\cG$ from $e^\prime\in\cE$ to
$e\in\cE$ is the tuple $\{e^\prime,e\}$ with
$\partial(e)=\partial(e^\prime)$. Both the combinatorial and the metric
length of a trivial walk are zero.

A walk $\bw=\{e^\prime,i_1,\ldots,i_N,e\}$ \emph{traverses} an internal
line $i\in\cI$ if $i_k=i$ for some $1\leq k \leq N$. It \emph{visits} the
vertex $v$ if either $v=\partial(e^\prime)$ or $v=\partial(e)$ or $v$ is
incident with at least one internal line traversed by the walk $\bw$.

The \emph{score} $\underline{n}(\bw)$ of a walk $\bw$ is the set
$\{n_i(\bw)\}_{i\in\cI}$ with $n_i(\bw)\geq 0$ being the number of times
the walk $\bw$ traverses the internal line $i\in\cI$. Any trivial walk has
the score $\underline{n}=\underline{0}:=\{0,\ldots,0\}$.

Let $\cW_{e,e^{\prime}} = \cW_{e,e^{\prime}}(\cG)$, $e,e^{\prime}\in\cE$ be
the set of all walks $\bw$ on $\cG$ from $e^{\prime}$ to $e$. In
particular, the set $\cW_{e,e^{\prime}}$ is infinite for all
$e,e^{\prime}\in\cE$ if $\cI\neq\emptyset$ and the graph $\cG$ is
connected. By reversing a walk $\bw$ from $e^\prime$ to $e$ into a walk
$\bw_{\mathrm{rev}}$ from $e$ to $e^\prime$ we obtain a natural one-to-one
correspondence between $\cW_{e,e^{\prime}}$ and $\cW_{e^{\prime},e}$.
Obviously, $|\bw|=|\bw_{\mathrm{rev}}|$ and
$\underline{n}(\bw)=\underline{n}(\bw_{\mathrm{rev}})$.

Let $\cS(v)\subseteq \cE\cup\cI$ denote the \emph{star graph} of the vertex
$v\in V$, i.e., the set of the edges adjacent to $v$. Also, by $\cS_{-}(v)$
(respectively $\cS_{+}(v)$) we denote the set of the edges for which $v$ is
the initial vertex (respectively terminal vertex). Obviously,
$\cS_{+}(v)\cap \cS_{-}(v)=\emptyset$ since $\cG$ does not contain
tadpoles by assumption.

To every $v\in V$ we associate an arbitrary $\deg(v)\times\deg(v)$ matrix
$M(v)$ with complex entries $[M(v)]_{j_1,j_2}$, where $j_1,j_2\in S(v)$ are
edges incident with the vertex $v$. The collection of such matrices for all
$v\in V$ will be denoted by $\cM=\{M(v)\}_{v\in V(\cG)}$.

Now to each non-trivial walk $\bw=\{e^\prime,i_1,\ldots,i_N,e\}$ from $e^\prime\in\cE$
to $e\in\cE$ on the graph $\cG$ we associate a weight $W(\bw)$ by
\begin{equation}\label{weight:2}
W(\bw)= \left[M(v_{|\bw_{\mathrm{comb}}|})\right]_{e,i_{|\bw|_{\mathrm{comb}}}}\cdot
\prod_{k=1}^{|\bw|_{\mathrm{comb}}-1}\left[M({v_{k}})\right]_{i_{k+1},
i_{k}}\,\cdot
\left[M(v_0)\right]_{i_1, e^{\prime}},
\end{equation}
where $v_0=\partial(e^\prime)$, $v_{|\bw|_{\mathrm{comb}}}=\partial(e)$,
$v_k$ with $k\in\{1,\ldots,|\bw|_{\mathrm{comb}}-1\}$ is the vertex incident
with the internal line $i_k$ as well as the internal line $i_{k+1}$. To a
trivial walk $\bw=\{e^\prime,e\}$ we associate the weight
\begin{equation}\label{weight:21}
W(\bw)=\left[M(\partial(e))\right]_{e,e^{\prime}}.
\end{equation}
\begin{definition}\label{def:gen:func}
The generating function of walks from $e^\prime\in\cE$ to $e\in\cE$ on the
graph $\cG$ associated with the collection $\cM=\{M(v)\}_{v\in V}$ is
defined as
\begin{equation}\label{statesum}
T_{e,e^{\prime}}(\beta) = \sum_{\bw\in\cW_{e,e^{\prime}}} W(\bw)
\e^{-\beta|\bw|} = \sum_{\bw\in\cW_{e,e^{\prime}}} W(\bw) \e^{-\beta\langle
\underline{n}(\bw), \underline{a}\rangle},
\end{equation}
where
\begin{equation*}
|\bw|=\langle \underline{n}(\bw), \underline{a}\rangle := \sum_{i\in\cI} n_i(\bw)
a_i.
\end{equation*}
\end{definition}

For given $\cM$ a walk $\bw$ is called \emph{relevant} if $W(\bw)\neq 0$.
The set of relevant walks from $e^{\prime}$ to $e$ is denoted by $\cW_{e,e^{\prime}}(\cM)$.

\begin{proposition}\label{prop:2.4}
There is $\beta_0>0$ such that the series \eqref{statesum} converges for
any $e,e^\prime\in\cE$ and all $\beta\in\C$ with $\Re\beta>\beta_0$.
Moreover,
\begin{equation}\label{tinfty}
\lim_{\Re\beta\uparrow\infty}T_{e,e^{\prime}}(\beta)=\begin{cases}
\left[M(\partial(e))\right]_{e,e^{\prime}}\quad & \text{if}\quad \partial(e)=\partial(e^{\prime}),\\
0\quad& \text{otherwise}.
\end{cases}
\end{equation}
\end{proposition}

Definition \ref{def:gen:func} suggests that we write $\cW_{e,e^{\prime}}$
as an infinite union of disjoint, non-empty sets by grouping together all
walks $\bw$ with the same score $\underline{n}(\bw)$,
\begin{equation*}
\cW_{e,e^{\prime}}(\underline{n})=\left\{{\bf w}\in
\cW_{e,e^{\prime}}\left.\right|\underline{n}({\bf w})=\underline{n}\right\}
\end{equation*}
such that
\begin{equation}\label{union}
\cW_{e,e^{\prime}}=\bigcup_{\underline{n}}\cW_{e,e^{\prime}}(\underline{n}).
\end{equation}
Note that these sets depend only on topology of the graph $\cG$ and are
independent of its metric properties. Also if $\bw\in
\cW_{e,e^{\prime}}(\underline{n})$ then $\bw_{\mathrm{rev}}\in
\cW_{e^{\prime},e}(\underline{n})$.
$\cW_{e,e^{\prime}}(\underline{0})=\emptyset$ if and only if $\partial(e) \neq
\partial(e^{\prime})$.

For the proof of Proposition \ref{prop:2.4} we need the following rather
obvious fact:

\begin{lemma}\label{4:lem:2}
The sets $\cW_{e,e^{\prime}}(\underline{n})$ are finite. Let
\begin{equation*}
|\underline{n}| =\sum_{i\in\cI}n_{i}
\end{equation*}
be the total number of internal lines traversed by any walk
$\bw\in\cW_{e,e^{\prime}}(\underline{n})$. The number of different walks in
$\cW_{e,e^{\prime}}(\underline{n})$ satisfies the bound
\begin{equation}
\label{nnumber} |\cW_{e,e^{\prime}}(\underline{n})|\leq
\frac{|\underline{n}|!} {\displaystyle\prod_{i\in\cI}n_{i}!}.
\end{equation}
\end{lemma}

Set
\begin{equation}\label{spar}
T_{e,e^{\prime}}(\underline{n}) =
\sum_{\bw\in\cW_{e,e^{\prime}}(\underline{n})}W(\bw)
\end{equation}
if $\cW_{e,e^{\prime}}(\underline{n})$ is nonempty and
$T_{e,e^{\prime}}(\underline{n})=0$ whenever
$\cW_{e,e^{\prime}}(\underline{n})=\emptyset$. Observe that
$T_{e,e^\prime}(\underline{n})$ does not depend on the metric properties of
the graph, i.e., is independent of the lengths of internal lines
$\underline{a}$.

For given $e,e^{\prime}\in\cE$ consider the set of scores of all walks from
$e^{\prime}$ to $e$,
\begin{equation}\label{enn}
\cN_{e,e^{\prime}}=\left\{\underline{n}\left.\right|
\text{there is a walk}\, \bw\in \cW_{e,e^{\prime}}(\underline{n})\right\}.
\end{equation}
Since $\underline{n}(\bw)=\underline{n}(\bw_{\mathrm{rev}})$, we have
$\cN_{e,e^{\prime}}=\cN_{e^{\prime},e}$.

With this notation we have the following equivalent representation of
\eqref{statesum}:
\begin{equation}\label{statesum1}
T_{e,e^\prime}(\beta) = \sum_{\underline{n}\in\cN_{e,e^\prime}}
T_{e,e^\prime}(\underline{n})\; \e^{-\beta \langle\underline{n},
\underline{a}\rangle} .
\end{equation}
Obviously, the series in \eqref{statesum} converges absolutely if and only
if the series in \eqref{statesum} does.

\begin{proof}[Proof of Proposition \ref{prop:2.4}]
Observe that
\begin{equation*}
\Big|\sum_{\bw\in\cW_{e,e^\prime}(\underline{n})} W(\bw)
\e^{-\beta\langle\underline{n}, \underline{a}\rangle}\Big| \leq
\sum_{\bw\in\cW_{e,e^\prime}(\underline{n})}\left(\max_{v\in
V}\|M(v)\|\right)^{|\underline{n}|+1} \e^{-|\underline{n}|\Re\beta\,
a_{\mathrm{min}}},
\end{equation*}
where
\begin{equation}
\label{amin}
0< a_{\mathrm{min}} := \min_{i\in\cI} a_i.
\end{equation}
{}From Lemma \ref{4:lem:2} and using the identity
\begin{equation}\label{identity}
\sum_{\substack{\underline{n}\in\N_0^{|\cI|}\\ |\underline{n}|=N}}
\:\frac{|\underline{n}|!}{\displaystyle\prod_{i\in\cI}n_{i}!} =
|\cI|^{N},\qquad N\in\N
\end{equation}
we, therefore, obtain
\begin{equation*}
\left|\sum_{\underline{n}\in\cN_{e,e^\prime}} \Big(\sum_{\bw\in
\cW_{e,e^\prime}(\underline{n})} W(\bw) \e^{-\beta \langle\underline{n},
\underline{a}\rangle} \Big)\right| \leq \sum_{N=0}^\infty \left(\max_{v\in
V}\|M(v)\|\right)^{N+1} \e^{-N \Re\beta\, a_{\mathrm{min}}} |\cI|^N.
\end{equation*}
This series converges for all $\beta\in\C$ with $\Re\beta>\beta_0$, where
\begin{equation}\label{beta:0}
\beta_0 > \frac{1}{a_{\mathrm{min}}}\left(\max_{v\in V}\log\|M(v)\| + \log
|\cI|\right).
\end{equation}
\end{proof}

We mention also the following simple result:

\begin{lemma}[Time Reversal Invariance]
If all matrices $M(v)$ are symmetric then so
is the matrix $T(\beta)$ with matrix elements $T_{e,e^\prime}(\beta)$ for
all large $\Re\beta>0$. If all $M(v)$ are self-adjoint, then so is
$T(\beta)$ for all large $\beta>0$.
\end{lemma}

\begin{definition}
The family of matrices $\cM$ is called \emph{combinatorial} if every matrix
entry of every matrix $M(v)$ equals either zero or one.
\end{definition}

If $\cM$ is combinatorial, the weight $W(\bw)$ of an arbitrary walk $\bw$
is either zero or one and we have the following simple result.

\begin{lemma}\label{comb}
If $\cM$ is combinatorial and $\cW_{e,e^{\prime}}(\cM)$ finite then
\begin{equation*}
T_{e,e^{\prime}}(0)= |\cW_{e,e^{\prime}}(\cM)|,
\end{equation*}
i.e., the number of relevant walks from $e^\prime\in\cE$ to $e\in\cE$.
\end{lemma}

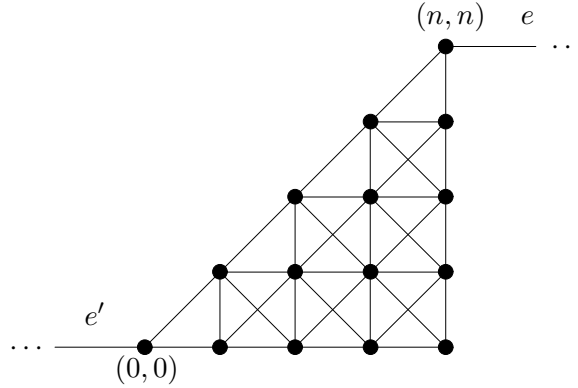
\begin{figure}
\setlength{\unitlength}{2mm}
\begin{center}
\begin{picture}(30,30)(-1,-1)
\matrixput(0,0)(5,0){5}(0,5){1}{\circle*{1}}
\matrixput(5,5)(5,0){4}(0,5){1}{\circle*{1}}
\matrixput(10,10)(5,0){3}(0,5){1}{\circle*{1}}
\matrixput(15,15)(5,0){2}(0,5){1}{\circle*{1}}
\matrixput(20,20)(5,0){1}(0,5){1}{\circle*{1}} \path(0,0)(20,0)
\path(0,0)(20,20) \path(20,20)(20,0) \path(5,5)(20,5) \path(10,10)(20,10)
\path(15,15)(20,15) \path(5,0)(5,5) \path(10,0)(10,10) \path(15,0)(15,15)
\path(5,5)(10,0) \path(10,10)(20,0) \path(15,15)(20,10) \path(10,5)(15,0)
\path(15,10)(20,5) \path(5,0)(20,15) \path(10,0)(20,10)
\path(15,0)(20,5)
\path(-6,0)(0,0)
\path(20,20)(26,20)
\put(-9,-0.5){$\cdots$}
\put(27,19.5){$\cdots$}
\put(-2,-2){$(0,0)$} \put(18,21.5){$(n,n)$}
\put(-4,1.5){$e^{\prime}$}
\put(25,21.5){$e$}
\end{picture}
\end{center}
\caption{\label{fig:1} The graph $\cG_n$ for $n=4$. The
metric lengths of all internal lines are assumed to be equal $1$. The
external lines are $e$ and $e^{\prime}$.}
\end{figure}

We now provide some examples, which relate our formulation to well known
combinatorial contexts. Viewing $\Z^2$ as a subset of $\R^2$, for an
arbitrary $n\in\N$ consider the set
\begin{equation*}
V_n = \left\{(x_1,x_2)\in\Z^2\big|\; 0\leq x_2\leq x_1\leq n\right\}.
\end{equation*}
We consider the non-compact graph $\cG_n=(V_n,\cI,\cE,\partial)$, where
$\cE=\{e,e^\prime\}$, $\partial(e^\prime)=(0,0)$ and $\partial(e)=(n,n)$,
and the vertices $v_1\in V_n$ and $v_2\in V_n$ are adjacent if and only if
the Euclidean distance between these vertices is not larger than $\sqrt{2}$,
$|v_1 - v_2|\leq \sqrt{2}$. Therefore, the set of internal lines $\cI$
consists of of all intervals joining the points of $V_n$ and having
Euclidean distance not greater than $\sqrt{2}$ (see Fig.~\ref{fig:1}). The
metric distance between two adjacent vertices will be assumed to be equal
$1$, that is, $a_i=1$ for all $i\in\cI$.

\begin{example}[The Catalan numbers]
The number
\begin{equation*}
C_{n-1}= \frac{1}{n+1}\begin{pmatrix}
2n \\ n
                      \end{pmatrix},\qquad n\in\N
\end{equation*}
is called the $(n-1)$-th Catalan number (see, e.g., \cite{West} and pp.~219
-- 229 in \cite{Stanley:2}). Set $\cK^{\textsl{Catalan}}=\{(1,0),(0,1)\}$.
For an arbitrary vertex $v\in V_n$ of the graph $\cG_n$ and arbitrary
$j\in\cE\cap\cI$ adjacent to the vertex $v$ we set
\begin{equation*}
\chi_v(j)=\begin{cases}
v^ \prime-v\in\Z^2 & \text{if}\; j\in\cI\; \text{and is adjacent to the vertex}\; v^ \prime,\\
(1,0)\in\Z^2 & \text{if}\; j=e,\\
(-1,0)\in\Z^2 & \text{if}\; j=e^ \prime.
          \end{cases}
\end{equation*}
Let $\cK^{\textsl{Catalan}}=\{(1,0),(0,1)\}\subset\Z^2$ and
\begin{equation*}
\left[M^{\textsl{Catalan}}(v)\right]_{j_1,j_2}=
\begin{cases}
1 & \text{if}\; \chi_v(j_1)\in\cK^{\textsl{Catalan}}\;\text{and}\; -\chi_v(j_2)\in\cK^{\textsl{Catalan}},\\
0 & \text{otherwise.}
\end{cases}
\end{equation*}

The set $\cW_{e,e^{\prime}}(\cM^{\textsl{Catalan}})$ is, obviously, finite.
Therefore, the generating function $T_{e,e^{\prime}}(\beta)$ is entire. For
given $n\in\N$ the number $T_{e,e^{\prime}}(0)$ is the $(n-1)$-th Catalan
number. The three other matrix elements $T_{e,e}(\beta)$,
$T_{e^\prime,e}(\beta)$, and $T_{e^\prime,e^\prime}(\beta)$ vanish
identically.
\end{example}

In the next example we continue with the same notation.

\begin{example}[The Schr\"{o}der numbers]
The Schr\"{o}der numbers (see, e.g., \cite{West} and p.~178 in \cite{Stanley:2})
can be defined by the recurrence relation
\begin{equation*}
S_n = S_{n-1} + \sum_{k=0}^{n-1} S_k S_{n-k-1}\quad\text{with}\quad S_0=1.
\end{equation*}
Let $\cK^{\textsl{Schr\"{o}der}}=\{(1,0),(0,1),(1,1)\} \supset
\cK^{\textsl{Catalan}}$ and
\begin{equation*}
\left[M^{\textsl{Schr\"{o}der}}(v)\right]_{j_1,j_2}=
\begin{cases}
1 & \text{if}\; \chi_v(j_1)\in\cK^{\textsl{Schr\"{o}der}}\;\text{and}\; -\chi_v(j_2)\in\cK^{\textsl{Schr\"{o}der}},\\
0 & \text{otherwise.}
\end{cases}
\end{equation*}

Obviously, $\cW_{e,e^{\prime}}(\cM^{\textsl{Catalan}})\subseteq
\cW_{e,e^{\prime}}(\cM^{\textsl{Schr\"{o}der}})$ is again a finite set. For
given $n\in\N$ the number $T_{e,e^{\prime}}(0)$ is now the $n$-th Schr\"{o}der
number. The three other matrix elements $T_{e,e}(\beta)$,
$T_{e^\prime,e}(\beta)$, and $T_{e^\prime,e^\prime}(\beta)$ vanish
identically.
\end{example}

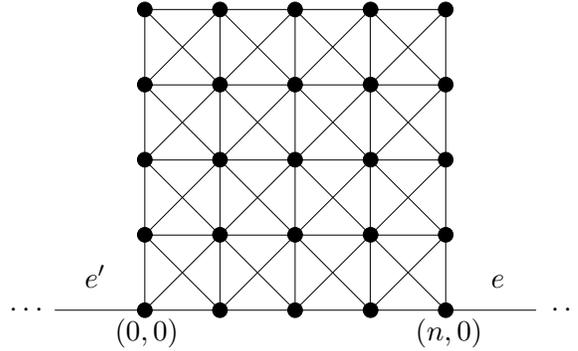
\begin{figure}
\setlength{\unitlength}{2mm}
\begin{center}
\begin{picture}(30,30)(-1,-1)
\matrixput(0,0)(5,0){5}(0,5){5}{\circle*{1}} \path(0,0)(20,20)
\path(0,0)(20,0) \path(0,0)(0,20) \path(20,0)(20,20) \path(0,20)(20,20)
\path(0,20)(20,0) \path(0,5)(20,5) \path(0,10)(20,10) \path(0,15)(20,15)
\path(5,0)(5,20) \path(10,0)(10,20) \path(15,0)(15,20) \path(5,0)(20,15)
\path(10,0)(20,10) \path(15,0)(20,5) \path(0,5)(15,20) \path(0,10)(10,20)
\path(0,15)(5,20) \path(0,5)(5,0) \path(0,10)(10,0) \path(0,15)(15,0)
\path(5,20)(20,5) \path(10,20)(20,10) \path(15,20)(20,15)
\put(-2,-2){$(0,0)$} \put(18,-2){$(n,0)$}
\path(-6,0)(0,0)
\path(20,0)(26,0)
\put(-9,-0.4){$\cdots$}
\put(27,-0.4){$\cdots$}
\put(-4,1.5){$e^{\prime}$}
\put(23,1.5){$e$}
\end{picture}
\end{center}
\caption{\label{fig:2} The graph $\cG_n^{+}$ for $n=4$.
The metric lengths of all internal lines are assumed to be equal
$1$. There are again two external lines $e$ and $e^{\prime}$. \hfill}
\end{figure}

For the next two examples consider the sets
\begin{equation*}
V_n^{+} = \left\{(x_1,x_2)\in\Z^2\big|\; 0\leq x_1\leq n,\; 0\leq x_2\leq
n\right\},\qquad n\in\N.
\end{equation*}
Let $\cG_n^{+} = (V_n^{+},\cI,\cE,\partial)$ be the non-compact graph with
$\cE=\{e, e^\prime\}$, $\partial(e^\prime)=(0,0)$ and $\partial(e)=(n,0)$,
and the vertices $v_1\in V$ and $v_2\in V$ are adjacent if and only if the
Euclidean distance between these vertices is not larger than $\sqrt{2}$,
$|v_1 - v_2|\leq \sqrt{2}$. Therefore, the set of internal lines $\cI$
consists of of all intervals joining the points of $V_n^ {+}$ and having
Euclidean length not greater than $\sqrt{2}$ (see Fig.~\ref{fig:2}).

\begin{example}[Dyck paths]\label{ex:2.9}
Let $\cK^{\textsl{Dyck}}=\{(1,1),(1,-1)\}$ and
\begin{equation*}
\left[M^{\textsl{Dyck}}(v)\right]_{j_1,j_2}=
\begin{cases}
1 & \text{if}\; \chi_v(j_1)\in\cK^{\textsl{Dyck}}\;\text{and}\; -\chi_v(j_2)\in\cK^{\textsl{Dyck}},\\
0 & \text{otherwise}
\end{cases}
\end{equation*}
if neither $j_1$ nor $j_2$ are external lines. We set
\begin{equation*}
\left[M^{\textsl{Dyck}}(v)\right]_{j_1,e^\prime}=
\begin{cases}
1 & \text{if}\; \chi_v(j_1)\in\cK^{\textsl{Dyck}},\\
0 & \text{otherwise}
\end{cases}
\end{equation*}
if $j_2=e^\prime$ and
\begin{equation*}
\left[M^{\textsl{Dyck}}(v)\right]_{e,j_2}=
\begin{cases}
1 & \text{if}\; -\chi_v(j_2)\in\cK^{\textsl{Dyck}},\\
0 & \text{otherwise}
\end{cases}
\end{equation*}
if $j_1=e$.

Obviously, $\cW_{e,e^\prime}(\cM^{\textsl{Dyck}})$ is a finite set.
Therefore, $T_{e,e^\prime}(\beta)$ is entire, $T_{e,e^\prime}(0)$ is the
number of Dyck paths on the graph $\cG_n^{+}$. A discussion of Dyck paths
can be found in \cite{Kra}.
\end{example}

\begin{example}[Motzkin numbers]
The non-compact graphs $\cG_n^{+}$ are the same as for Dyck paths in Example
\ref{ex:2.9}. The Motzkin numbers (see, e.g., \cite{Aigner} and Problem 6.37
in \cite{Stanley:2}) can be defined by the recurrence relation
\begin{equation*}
M_n = M_{n-1} + \sum_{k=0}^{n-2} M_k M_{n-k-2}\quad\text{with}\quad
M_0=M_1=1.
\end{equation*}
Set $\cK^{\textsl{Motzkin}}=\{(1,1), (1,-1), (1,0)\}\supset
\cK^{\textsl{Dyck}}$ and \begin{equation*}
\left[M^{\textsl{Motzkin}}(v)\right]_{j_1,j_2}=
\begin{cases}
1 & \text{if}\; \chi_v(j_1)\in\cK^{\textsl{Motzkin}}\;\text{and}\; -\chi_v(j_2)\in\cK^{\textsl{Motzkin}},\\
0 & \text{otherwise.}
\end{cases}
\end{equation*}

Again the set $\cW_{e,e^\prime}(\cM^{\textsl{Motzkin}})$ is finite and,
therefore, $T_{e,e^\prime}(\beta)$ is entire. For given $n\in\N$ the number
$T_{e,e^{\prime}}(0)$ is the $n$-th Motzkin number.
\end{example}

\section{Laplace Operators on Graphs}\label{sec:operator}

In this section we will recall the theory of Laplace operators on a metric
graph $\cG$ and the resulting scattering theory (see \cite{KS1},
\cite{KS2}, \cite{KS3}, \cite{KS4}, \cite{KS5} for further details).

Given a finite non-compact graph $\cG=(V,\cI,\cE,\partial)$ with a metric
structure $\underline{a}=\{a_i\}_{i\in\cI}$ consider the Hilbert space
\begin{equation}\label{hilbert}
\cH\equiv\cH(\cE,\cI,\underline{a})=\cH_{\cE}\oplus\cH_{\cI},\qquad
\cH_{\cE}=\bigoplus_{e\in\cE}\cH_{e},\qquad
\cH_{\cI}=\bigoplus_{i\in\cI}\cH_{i},
\end{equation}
where $\cH_{e}=L^2([0,\infty))$ for all $e\in\cE$ and
$\cH_{i}=L^2([0,a_{i}])$ for all $i\in\cI$. By $\cD_j$ with
$j\in\cE\cup\cI$ denote the set of all $\psi_j\in\cH_j$ such that
$\psi_j(x)$ and its derivative $\psi^\prime_j(x)$ are absolutely continuous
and $\psi^{\prime\prime}_j(x)$ is square integrable. Let $\cD_j^0$ denote
the set of those elements $\psi_j$ of $\cD_j$ which satisfy
\begin{equation*}
\psi_j(0)=\psi_j(a_j)=
\psi_j^{\prime}(0)=\psi_j^{\prime}(a_j)=0\quad\text{for}\quad j\in\cI
\end{equation*}
and
\begin{equation*}
\psi_j(0) = \psi_j^{\prime}(0)=0\quad\text{for}\quad j\in\cE.
\end{equation*}
Let $\Delta^0$ be the differential operator
\begin{equation}\label{Delta:0}
\left(\Delta^0\psi\right)_j (x) = \frac{d^2}{dx^2} \psi_j(x),\qquad
j\in\cI\cup\cE
\end{equation}
with $\psi=\{\psi_j\}_{j\in\cI\cup\cE}$ in the domain
\begin{equation*}
\cD^0=\bigoplus_{j\in\cE\cup\cI} \cD_j^0 \subset\cH.
\end{equation*}
It is straightforward to verify that $\Delta^0$ is a closed symmetric
operator with deficiency indices equal to $|\cE|+2|\cI|$.

We introduce an auxiliary finite-dimensional Hilbert space
\begin{equation}\label{K:def}
\cK\equiv\cK(\cE,\cI)=\cK_{\cE}\oplus\cK_{\cI}^{(-)}\oplus\cK_{\cI}^{(+)}
\end{equation}
with $\cK_{\cE}\cong\C^{|\cE|}$ and $\cK_{\cI}^{(\pm)}\cong\C^{|\cI|}$. The
subspaces $\cK_{\cI}^{(-)}$ we associate with initial vertices of the
internal lines $i\in\cI$, the subspaces $\cK_{\cI}^{(+)}$ with the terminal
vertices. Let ${}^d\cK$ denote the ``double'' of $\cK$, that is,
${}^d\cK=\cK\oplus\cK$.

For any $\displaystyle\psi\in\cD:=\bigoplus_{j\in\cE\cup\cI} \cD_j$ we set
\begin{equation}\label{lin1}
[\psi]:=\underline{\psi}\oplus \underline{\psi}^\prime\in{}^d\cK,
\end{equation}
with
\begin{equation}\label{lin1:add}
\underline{\psi} = \begin{pmatrix} \{\psi_e(0)\}_{e\in\cE} \\
                                   \{\psi_i(0)\}_{i\in\cI} \\
                                   \{\psi_i(a_i)\}_{i\in\cI} \\
                                     \end{pmatrix}\in\cK,\qquad
\underline{\psi}^\prime = \begin{pmatrix} \{\psi_e^\prime(0)\}_{e\in\cE} \\
                                   \{\psi_i^\prime(0)\}_{i\in\cI} \\
                                   \{-\psi_i^\prime(a_i)\}_{i\in\cI} \\
                                     \end{pmatrix}\in\cK.
\end{equation}
Here the vector notation is used with respect to the orthogonal
decomposition \eqref{K:def}.

To define the Laplace operator on the graph $\cG$ consider the family
$\psi=\{\psi_{j}\}_{j\in\cE\cup\cI}$ of complex valued functions defined on
$[0,\infty)$ if $j\in\cE$ and on $[0,a_{i}]$ if $j\in\cI$. Formally the
(self-adjoint) Laplace operator is defined as
\begin{equation}\label{diff:expression}
\left(\Delta(A,B,\underline{a})\psi\right)_j (x) = \frac{d^2}{dx^2}
\psi_j(x),\qquad j\in\cI\cup\cE
\end{equation}
with the boundary conditions
\begin{equation}\label{Randbedingungen}
A\underline{\psi} + B\underline{\psi}' = 0.
\end{equation}
By definition $A$ and $B$ are any complex
$(|\cE|+2|\cI|)\times(|\cE|+2|\cI|)$ matrices
such that
\begin{equation}\label{abcond}
\begin{split}
\text{(i)} &\quad \text{the matrix $(A, B)$ has maximal rank},\\
\text{(ii)} &\quad \text{the matrix $A B^\dagger$ is self-adjoint}.
\end{split}
\end{equation}
Here and in what follows $(A,B)$ will denote the $(|\cE|+2|\cI|)\times
2(|\cE|+2|\cI|)$ matrix, where $A$ and $B$ are put next to each other.

The scattering matrix $S(\sk)=S(\sk;A,B,\underline{a})$ associated to
$\Delta(A,B,\underline{a})$ has the following interpretation in terms of
the solutions to the Schr\"{o}dinger equation (see \cite{KS1} and \cite{KS4}).
Consider the solutions $\psi^{k}(\sk)\;(k\in \cE)$ of the stationary
Schr\"{o}dinger equation for $-\Delta(A,B,\underline{a})$ at energy $\sk^2>0$,
\begin{equation*}
-\Delta(A,B,\underline{a})\psi^{k}(\sk)=\sk^2\psi^{k}(\sk)
\end{equation*}
of the form
\begin{equation}
\label{10} \psi^{k}_{j}(x;\sk)=\begin{cases} S(\sk)_{jk} \e^{\i\sk x}
&\text{for}
                                          \;j\in\cE, j\neq k\\
  \e^{-\i\sk x}+S(\sk)_{kk} \e^{\i\sk x} & \text{for}\;j\in\cE, j=k\\
                                  \alpha(\sk)_{jk} \e^{\i\sk x}+
        \beta(\sk)_{jk} \e^{-\i\sk x} & \text{for}\; j\in\cI. \end{cases}
\end{equation}
Thus, the number $S(\sk)_{jk}$ for $j\neq k$ is the transmission amplitude
from channel $k\in\cE$ to channel $j\in\cE$ and $S(\sk)_{kk}$ is the
reflection amplitude in channel $k\in\cE$. Their absolute squares may be
interpreted as transmission and reflection probabilities, respectively. The
``interior'' amplitudes
\begin{equation*}
\alpha(\sk)_{jk}=\alpha(\sk;A,B,\underline{a})_{jk},\qquad
\beta(\sk)_{jk}=\beta(\sk;A,B,\underline{a})_{jk}
\end{equation*}
are also of interest, since they describe how an incoming wave moves through
a graph before it is scattered into an outgoing channel.

The condition for the $\psi^{k}(E)\;(k\in \cE)$ to satisfy the boundary
conditions \eqref{Randbedingungen} immediately leads to the following
solution for the scattering matrix $S(\sk):\;
\cK_{\cE}\rightarrow\cK_{\cE}$ and the operators $\alpha(\sk)$ and
$\beta(\sk)$ acting from $\cK_{\cE}$ to $\cK_{\cI}$. Indeed, by combining
these operators into a map $\cK_{\cE}$ to
$\cK=\cK_{\cE}\oplus\cK_{\cI}^{(-)}\oplus\cK_{\cI}^{(+)}$ we obtain the linear equation
\begin{equation}\label{11}
Z(\sk;A,B,\underline{a})\begin{pmatrix} S(\sk)\\
                      \alpha(\sk)\\
                    \beta(\sk)\end{pmatrix}
                      =-(A-\i\sk B)
              \begin{pmatrix} \1\\
                               0\\
                               0 \end{pmatrix}
\end{equation}
with
\begin{equation}\label{Z:def}
Z(\sk;A,B,\underline{a})= A X(\sk;\underline{a})+\i\sk B
Y(\sk;\underline{a}),
\end{equation}
where
\begin{equation}\label{zet}
X(\sk;\underline{a}) = \begin{pmatrix}\1&0&0\\
                                  0&\1&\1\\
               0&\e^{\i\sk\underline{a}}&\e^{-\i\sk\underline{a}}
               \end{pmatrix} ,\qquad
Y(\sk;\underline{a}) = \begin{pmatrix}\1&0&0\\
                                  0&\1&-\1\\
               0&-\e^{\i\sk\underline{a}}&\e^{-\i\sk\underline{a}}
               \end{pmatrix}.
\end{equation}
The diagonal $|\cI|\times |\cI|$ matrices $\e^{\pm \i\sk\underline{a}}$ are
given by
\begin{equation}
\label{diag}
 \e^{\pm \i\sk\underline{a}}_{jk}=\delta_{jk}\e^{\pm \i\sk a_{j}}\;
                       \text{for}\; j,k\in\;\cI.
\end{equation}

\begin{theorem}[= Theorem 3.2 in \cite{KS1}]\label{3.2inKS1}
For any $\sk\in\R$
\begin{equation*}
\Ran\; (A-\i\sk B) \begin{pmatrix} \1 \\ 0 \\ 0 \end{pmatrix}\subset \Ran
Z(\sk;A,B,\underline{a}).
\end{equation*}
Thus, equation \eqref{11} has a solution even if $\det
Z(\sk;A,B,\underline{a})=0$ for some $\sk\in\R$. This solution defines the
scattering matrix uniquely. Moreover,
\begin{equation}\label{S-matrix}
S(\sk)=-\begin{pmatrix} \1 & 0 & 0 \end{pmatrix}
Z(\sk;A,B,\underline{a})^{-1} P_{\Ker Z(\sk;A,B,\underline{a})}^\perp
(A-\i\sk B) \begin{pmatrix} \1 \\ 0 \\ 0 \end{pmatrix}
\end{equation}
is unitary for all $k\in\R\setminus\{0\}$.
\end{theorem}

In the case with no internal lines $(\cI=\emptyset)$ the relation
\eqref{S-matrix} for the scattering matrix simplifies to
\begin{equation}\label{svertex}
S(\sk;A,B)=-\left(A+\i\sk B\right)^{-1}\left(A-\i\sk
B\right).
\end{equation}

\begin{proposition}\label{proposition:2.3}
If $\det(A + \i\sk B)=0$ for some $\sk\in\C$, then $\sk=\i\kappa$ with
$\kappa\in\R$. For any sufficiently large $\rho>0$ there is a constant
$C_\rho>0$ such that
\begin{equation}\label{absch}
\|(A+\i\sk B)^{-1}\| \leq C_\rho (1+|\sk|)^{-1}
\end{equation}
for all $\sk\in\C$ with $|\sk|>\rho$.
\end{proposition}

\begin{proof}
Assume that $\det(A + \i\sk B)=0$ for some $\sk\in\C$ with $\Re \sk\neq 0$.
Then also
\begin{equation*}
\det(A^\dagger - \i\overline{\sk} B^\dagger)= \overline{\det(A + \i\sk
B)}=0.
\end{equation*}
Therefore, there is a $\chi\neq 0$ such that
\begin{equation}\label{neu:1}
(A^\dagger - \i\overline{\sk} B^\dagger)\chi=0.
\end{equation}
In particular, we have $(B A^\dagger - \i\overline{\sk} B
B^\dagger)\chi=0$. Therefore, since $B A^\dagger$ is self-adjoint, we get
\begin{equation*}
\begin{split}
&\langle\chi, BA^\dagger \chi\rangle = \langle\chi, BB^\dagger
\chi\rangle\, \Im\sk ,\\
&  \langle\chi, BB^\dagger \chi\rangle\, \Re\sk =0.
\end{split}
\end{equation*}
The second equality implies that $\chi\in\Ker B^\dagger$. Then, by
\eqref{neu:1}, $\chi\in\Ker A^\dagger$. Since the matrix $(A,B)$ is of
maximal rank, we have $\Ker A^\dagger\cap \Ker B^\dagger=\{0\}$. Thus,
$\chi=0$ which contradicts the assumption and, hence, $\Re \sk=0$.

Since $\det(A+\i\sk B)$ is a polynomial in $\sk$, it has a finite number of
zeroes. Take an arbitrary $\rho>0$ such that all its zeroes lie in the disk
$|\sk|<\rho$. Using the matrix inverse formula we represent any element of
$(A+\i\sk B)^{-1}$ as a quotient of two polynomials of degrees $|\cE| +
2|\cI|-1$ and $|\cE| + 2|\cI|$, respectively. In turn, this implies the
estimate \eqref{absch}.
\end{proof}

\begin{theorem}\label{thm:cont}
The scattering matrix $S(\sk)=S(\sk;A,B,\underline{a})$ is a meromorphic
function in the complex $\sk$-plane. In upper half-plane $\Im\sk>0$ it has at most a finite number
of poles which are located on the imaginary semiaxis $\Re\sk=0$.
Outside these poles the scattering matrix is holomorphic for all $\Im\sk>0$
and determined by the relation
\begin{equation}\label{S-matrix:analyt}
S(\sk)=-\begin{pmatrix} \1 & 0 & 0 \end{pmatrix}
Z(\sk;A,B,\underline{a})^{-1} (A-\i\sk B) \begin{pmatrix} \1 \\ 0 \\ 0
\end{pmatrix}.
\end{equation}
\end{theorem}

\begin{proof}
Assume that $\det Z(\sk;A,B,\underline{a})=0$ for some $k\in\C$ with
$\Im\sk>0$ and $\Re\sk\neq 0$. This implies that the homogeneous equation
\begin{equation*}
Z(\sk;A,B,\underline{a})\begin{pmatrix} s \\ \alpha \\ \beta \end{pmatrix} = 0
\end{equation*}
has a nontrivial solution with $s\in\cK_{\cE}$ and
$\alpha,\beta\in\C^{|\cI|}$. Consider the function
$\psi(x)=\{\psi_j(x)\}_{j\in\cI\cup\cE}$ defined by
\begin{equation*}
\psi_{j}(x)=\begin{cases} s_{j} \e^{\i\sk x} &\text{for}
                                          \;j\in\cE, \\
                                  \alpha_{j} \e^{\i\sk x}+
        \beta_{j} \e^{-\i\sk x} & \text{for}\; j\in\cI. \end{cases}
\end{equation*}
Obviously, $\psi(x)$ satisfies the boundary conditions
\eqref{Randbedingungen}. Moreover, $\psi\in L^2(\cG)$ since $\Im\sk>0$.
Hence, $\sk^2\in\C$ with $\Im\sk^2\neq 0$ is an eigenvalue of the operator
$\Delta(A,B,\underline{a})$ which contradicts the self-adjointness of
$\Delta(A,B,\underline{a})$.

Since $\det Z(\sk;A,B,\underline{a})$ is an entire function in $\sk$ which
does not vanish identically, from \eqref{11} it follows that the scattering
matrix $S(\sk)$ is a meromorphic function in the complex $\sk$-plane. To
prove that the scattering matrix $S(\sk)$ has at most a finite number of
poles on the imaginary semiaxis $\{\sk\in\C|\;\Re\sk=0,\; \Im\sk>0\}$ it
suffices to show that the determinant $\det Z(\sk;A,B,\underline{a})$ does
not vanish for all sufficiently large $\Im\sk>0$. To see this we set
$\sk=\i\kappa$ with $\kappa>0$ and assume there is an unbounded
non-decreasing sequence $\{\kappa_k\}_{k\in\N}$ such that
\begin{equation*}
\det Z(\i\kappa_k;A,B,\underline{a}) = 0\qquad\text{for all}\quad k\in\N.
\end{equation*}
Therefore, there is a sequence $\{\chi_k\}_{k\in\N}$ of normalized elements
$\chi_k\in\cK$ such that
\begin{equation*}
X(\i\kappa_k;\underline{a})^\dagger A^\dagger\chi_k=\kappa_k
Y(\i\kappa_k;\underline{a})^\dagger B^\dagger\chi_k.
\end{equation*}
It is straightforward to verify that $X(\i\kappa_k;\underline{a})$ is
invertible and
\begin{equation*}
R_k:=\left(X(\i\kappa_k;\underline{a})^\dagger\right)^{-1}Y(\i\kappa_k;\underline{a})^\dagger
=\begin{pmatrix}
\1 & 0 & 0 \\ 0 & \coth(\kappa \underline{a}) &  - [\sinh(\kappa \underline{a})]^{-1}\\
0 &  -[\sinh(\kappa \underline{a})]^{-1} & \coth(\kappa \underline{a})
 \end{pmatrix}
\end{equation*}
with a notation analogous to \eqref{diag}. Thus,
\begin{equation}\label{bbbb}
(A^\dagger-\kappa_k B^\dagger)\chi_k = \kappa_k(R_k-\1)B^\dagger\chi_k
\end{equation}
for all $k\in\N$. Observe that $\|R_k-\1\| = O(e^{-c\kappa_k})$ for some
$c>0$ as $k\rightarrow\infty$. By Proposition \ref{proposition:2.3} the
operator $A^\dagger-\kappa B^\dagger$ is invertible for all sufficiently
large $\kappa$. Moreover, $\|(A^\dagger-\kappa B^\dagger)^{-1}\|\leq C$
with $C>0$ for all sufficiently large $\kappa$. Thus, equation \eqref{bbbb}
implies that $\chi_k\rightarrow 0$ which contradicts the assumption
$\|\chi_k\|=1$.
\end{proof}

In the lower half-plane $\Im\sk <0$ the scattering matrix may have poles
with $\Re\sk\neq 0$ (see, e.g., Example 3.2 in \cite{KS1}). These poles
correspond to resonances.

The notion of local boundary conditions has been introduced in our article
\cite{KS1} and is discussed in more details in \cite{KS4} and \cite{KS5}.
Local boundary conditions couple only those boundary values of $\psi$ and
of its derivative $\psi^\prime$ which belong to the same vertex. The precise
definition is as follows.

With respect to the orthogonal decomposition $\cK =
\cK_{\cE}\oplus\cK_{\cI}^{(-)}\oplus\cK_{\cI}^{(+)}$ any element $z$ of
$\cK$ can be represented as a vector
\begin{equation}\label{elements}
z=\begin{pmatrix}\{z_e\}_{e\in\cE}\\ \{z^{(-)}_i\}_{i\in\cI}\\
\{z^{(+)}_i\}_{i\in\cI}\end{pmatrix}.
\end{equation}
Consider the orthogonal decomposition
\begin{equation*}
\cK = \bigoplus_{v\in V} \cL_{v}
\end{equation*}
with $\cL_{v}$ being the linear subspace of dimension $\deg(v)$ spanned by
those elements \eqref{elements} of $\cK$ which satisfy
\begin{equation}
\label{decomp}
\begin{split}
z_e=0 &\quad \text{if}\quad e\in \cE\quad\text{is not incident with the vertex}\quad v,\\
z^{(-)}_i=0 &\quad \text{if}\quad v\quad\text{is not an initial vertex of}\quad i\in \cI,\\
z^{(+)}_i=0 &\quad \text{if}\quad v\quad\text{is not a terminal vertex
of}\quad i\in \cI.
\end{split}
\end{equation}

Set ${^d}\cL_v:=\cL_v\oplus\cL_v\cong\C^{2\deg(v)}$. By the First Theorem of
Graph Theory we have
\begin{equation*}
\sum_{v\in V(\cG)} \deg(v) = |\cE| + 2|\cI|
\end{equation*}
such that
\begin{equation*}
\bigoplus_{v\in V(\cG)} {^d}\cL_v = {}^d\cK.
\end{equation*}

\begin{definition}\label{def:local}
Given the graph $\cG=\cG(V,\cI,\cE,\partial)$, the boundary conditions
$(A,B)$ satisfying \eqref{abcond} are called \emph{local} on $\cG$ if and
only if there is an invertible map $C:\, \cK\rightarrow\cK$ and linear
transformations $A(v)$ and $B(v)$ in $\cL_{v}$ such that the direct sum
decompositions
\begin{equation}\label{permut}
CA= \bigoplus_{v\in V} A(v)\quad \text{and}\quad CB= \bigoplus_{v\in V} B(v)
\end{equation}
hold simultaneously. Otherwise the boundary conditions are called \emph{non-local}.
\end{definition}

For instance, for a single-vertex graph any boundary conditions are local. The boundary conditions
considered in Example 3.4 of \cite{KS5} are non-local.

\section{Combinatorial Fourier Expansion of the Scattering Matrix}\label{sec:harmony}

In this section we will perform a harmonic analysis of the scattering matrix
with respect to the lengths
$\underline{a}=\{a_i\}_{i\in\cI}\in(\R_+)^{|\cI|}$ of the internal lines of
the graph $\cG$. The main results of this section are presented in Theorems
\ref{thm:main:harmony} and \ref{verloren}. In Theorem
\ref{thm:main:harmony} the absolute convergence of the Fourier series for
the scattering matrix is proved. Theorem \ref{verloren} expresses its
Fourier coefficients as sums over the walks on the graph. Combining these
two results proves the combinatorial Fourier expansion formula
\eqref{sseries}.

Throughout the whole section we will assume that the (topological) graph
$\cG$ as well as the boundary conditions $(A,B)$ are fixed. To carry out
the analysis we will now treat $\underline{a}$ as a parameter which may
belong to $\R^{|\cI|}$ or even $\C^{|\cI|}$.

We start with the following simple but important observation.

\begin{lemma}\label{periodic}
For arbitrary $\sk>0$ the scattering matrix $S(\sk;A,B,\underline{a})$ is
uniquely defined as a solution of \eqref{11} for all
$\underline{a}\in\R^{|\cI|}$. Moreover, the scattering matrix is periodic
with respect to $\underline{a}$,
\begin{equation*}
S\left(\sk;A,B,\underline{a}+\frac{2\pi}{\sk}\underline{\ell}\right)=S(\sk;A,B,\underline{a})
\end{equation*}
for arbitrary $\underline{\ell}\in\Z^{|\cI|}$.
\end{lemma}

\begin{proof}
It suffices to consider those $\underline{a}\in\R^{|\cI|}$ for which $\det
Z(\sk; A, B, \underline{a})=0$, since the claim is obvious when the
determinant is non-vanishing. For $\underline{a}\in(\R_+)^{|\cI|}$ the fact
that $S(\sk; A,B,\underline{a})$ is uniquely defined as a solution of
\eqref{11} is guaranteed by Theorem \ref{3.2inKS1}. The case of arbitrary
$\underline{a}\in\R^{|\cI|}$ can be treated exactly in the same way (see
the proof of Theorem 3.2 in \cite{KS1}).

The periodicity follows immediately from \eqref{11} and the fact that the
matrices $X(\sk;\underline{a})$ and $Y(\sk;\underline{a})$ in \eqref{zet}
are $\frac{2\pi}{\sk}\Z^{|\cI|}$-periodic.
\end{proof}

Lemma \ref{periodic} suggests to consider a Fourier expansion of the
scattering matrix. The following theorem ensures the absolute convergence
of the corresponding Fourier series.

\begin{theorem}\label{thm:main:harmony}
Let $\sk>0$ be arbitrary. For all $\underline{a}\in\R^{|\cI|}$ the Fourier
expansion of the scattering matrix
\begin{equation}\label{Fourier:exp}
S(\sk;A,B,\underline{a})=\sum_{\underline{n}\in\Z^{|\cI|}}
\widehat{S}_{\underline{n}}(\sk;A,B)\,
\e^{\i\sk\langle\underline{n},\underline{a}\rangle}
\end{equation}
with
\begin{equation}\label{fourier:coef}
\widehat{S}_{\underline{n}}(\sk;A,B) =
\left(\frac{\sk}{2\pi}\right)^{\mid\cI\mid}\;
\int\limits_{[0,2\pi/\sk]^{|\cI|}}d
\underline{a}\:S(\sk;A,B,\underline{a})\:
\e^{-\i\sk\langle\underline{n},\underline{a}\rangle}
\end{equation}
converges absolutely and uniformly on compact subsets of $\R^{|\cI|}$. The
Fourier coefficients \eqref{fourier:coef} vanish for all
$\underline{n}=\{n_i\}_{i\in\cI}\in\Z^{|\cI|}$ with $n_i<0$ for at
least one $i\in\cI$.
\end{theorem}

For the proof we need a couple of auxiliary results. Set
\begin{equation*}
\cA=\left\{\underline{a}=\{a_i\}_{i\in\cI}\big|\: \Re a_i\in\R,\:\Im a_i
>0\right\}\subset\C^{|\cI|}.
\end{equation*}

\begin{lemma}\label{lem:det:neq:0}
For any $\sk>0$ the determinant $\det Z(\sk;A,B,\underline{a})$ has no
zeroes for all $\underline{a}\in\cA$.
\end{lemma}

\begin{proof}
Assume there is $\underline{a}\in \cA$ such that $\det
Z(\sk;A,B,\underline{a})=0$. Then there are $s\in\C^{|\cE|}$ and $\alpha,
\beta\in\C^{|\cI|}$ such that
\begin{equation*}
Z(\sk;A,B,\underline{a})\begin{pmatrix} s \\ \alpha \\ \beta\end{pmatrix} =
0.
\end{equation*}
Equivalently this gives
\begin{equation*}
(A+\i\sk B)\begin{pmatrix} s \\ \alpha \\
\e^{-\i\sk\underline{a}}\beta \end{pmatrix} + (A-\i\sk B)\begin{pmatrix} 0 \\ \beta \\
\e^{\i\sk\underline{a}}\alpha \end{pmatrix}=0.
\end{equation*}
The operator $(A+\i\sk B)^{-1}(A-\i\sk B)$ is unitary for all $\sk>0$ (see
the proof of Theorem 2.1 in \cite{KS1} ). Since unitary transformations
preserve the canonical Hilbert norm on $\C^{|\cE|+2|\cI|}$, we have
\begin{equation*}
\|s\|^2 + \sum_{i\in\cI} |\alpha_i|^2 (1-\e^{-2\sk \Im a_i}) +
\sum_{i\in\cI} |\beta_i|^2 (\e^{2\sk \Im a_i}-1) = 0,
\end{equation*}
which implies $s=0$ and $\alpha=\beta=0$.
\end{proof}

\begin{proposition}\label{thm:6.1:harmo}
Let $\sk>0$ be arbitrary. For all
\begin{equation*}
\underline{a}\in\mathrm{clos}(\cA):=\left\{\underline{a}=\{a_i\}_{i\in\cI}\big|\:
\Re a_i\in\R,\: \Im a_i \geq 0\right\}
\end{equation*}
the scattering matrix $S(\sk;A,B,\underline{a})$ is uniquely defined as a
solution of \eqref{11} and satisfies the bound
\begin{equation}\label{S:leq:1}
\|S(\sk;A,B,\underline{a})\| \leq 1.
\end{equation}
Moreover, it is a rational function of $\underline{t}=\{t_i\}_{i\in\cI}$
with $t_i:=\e^{\i\sk a_i}$, i.e.\ a quotient of $\cB(\cK_{\cE})$-valued
polynomials in the variables $t_i$. Thus, for all
$\underline{a}\in\mathrm{clos}(\cA)$ the scattering matrix is
$\frac{2\pi}{\sk}\Z^{|\cI|}$-periodic,
\begin{equation*}
S\left(\sk;A,B,\underline{a}+\frac{2\pi}{\sk}\underline{\ell}\right)=S(\sk,A,B,\underline{a}),\qquad
\underline{\ell}\in\Z^{|\cI|}.
\end{equation*}
\end{proposition}

\begin{proof}
By Lemma \ref{lem:det:neq:0} equation \eqref{11} has a unique solution for
all $\underline{a}\in\cA$. Equations \eqref{Z:def} and \eqref{zet} imply
that $Z(\sk;A,B,\underline{a})$ is a polynomial function of the components
of $\underline{t}$. Obviously, $Z(\sk;A,B,\underline{a})^{-1}$ is also a
rational function of $\underline{t}$. Thus, by \eqref{11} the scattering
matrix $S(\sk;A,B,\underline{a})$ is a rational function of
$\underline{t}$. Thus, it is $\frac{2\pi}{\sk}\Z^{|\cI|}$-periodic.

Using \eqref{11} it is easy to check that this solution satisfies the
relation
\begin{equation*}
\begin{pmatrix} S(\sk;A,B,\underline{a})\\ \alpha(\sk;A,B,\underline{a})\\
\e^{-\i\sk\underline{a}}\beta(\sk;A,B,\underline{a}) \end{pmatrix} =
-(A+\i\sk B)^{-1} (A-\i\sk B) \begin{pmatrix} \1\\ \beta(\sk;A,B,\underline{a})\\
\e^{\i\sk\underline{a}}\alpha(\sk;A,B,\underline{a}) \end{pmatrix}.
\end{equation*}
Since $(A+\i\sk B)^{-1} (A-\i\sk B)$ is unitary we obtain
\begin{equation}\label{unitaer}
\begin{split}
& S(\sk;A,B,\underline{a})^\dagger S(\sk;A,B,\underline{a}) +
\alpha(\sk;A,B,\underline{a})^\dagger (\1-\e^{-2\sk \Im\underline{a}})
\alpha(\sk;A,B,\underline{a})\\
& \qquad + \beta(\sk;A,B,\underline{a})^\dagger (\e^{2\sk
\Im\underline{a}}-\1) \beta(\sk;A,B,\underline{a})=\1,
\end{split}
\end{equation}
where $\Im\underline{a}=\{\Im a_i\}_{i\in\cI}$. {}From \eqref{unitaer} it
follows immediately that
\begin{equation*}
0\leq S(\sk;A,B,\underline{a})^\dagger S(\sk;A,B,\underline{a}) \leq \1
\end{equation*}
in the operator sense. This proves the bound \eqref{S:leq:1} for all
$\underline{a}\in\cA$. Recalling Lemma \ref{periodic} completes the proof.
\end{proof}

A priori it is not clear whether the boundary values of the scattering
matrix $S(\sk; A,B,\underline{a})$ with $\underline{a}\in\cA$ coincide
with those given by equation \eqref{S-matrix} for all
$\underline{a}\in\R^{|\cI|}$. The following lemma shows the
``non-tangential continuity'' of the scattering matrix $S(\sk;A, B,
\underline{a})$ with respect to $\underline{a}\in\mathrm{clos}(\cA)$.

\begin{lemma}\label{Grenzwert}
Let $\underline{a}\in\R^{|\cI|}$ and $\sk>0$ be arbitrary. For any sequence
$\{\underline{a}_j\}_{j\in\N}$, $\underline{a}_j\in\cA$ converging to
$\underline{a}\in\R^{|\cI|}$ the relation
\begin{equation}\label{S:convergence}
\lim_{j\rightarrow\infty} S(\sk;A, B, \underline{a}_j) = S(\sk;A, B,
\underline{a})
\end{equation}
holds.
\end{lemma}

For the proof we need the following elementary result.

\begin{lemma}\label{Grenzwert:Hilfe}
Let $T_n$ be a sequence of invertible operators on the finite-dimensional
Hilbert space $\fH$ converging to the operator $T$. Then
\begin{equation*}
\lim_{n\rightarrow\infty} P_{\Ker T}^\perp T_n^{-1} P_{\Ker
T^\dagger}^\perp = P_{\Ker T}^\perp T^{-1} P_{\Ker T^\dagger}^\perp,
\end{equation*}
where $P_{\fL}^\perp$ denotes the orthogonal projection onto orthogonal
complement in $\fH$ of the subspace $\fL\subset\fH$.
\end{lemma}

\begin{proof}
Consider the operators $T_n$ and $T$ as maps from $(\Ker T)^\perp$ to
$(\Ker T^\dagger)^\perp$. Since these maps are invertible, the claim
follows from the obvious relation
\begin{equation*}
T_n^{-1} = T^{-1}\left[I+T^{-1}(T_n-T)\right]^{-1}.
\end{equation*}
\end{proof}

\begin{proof}[Proof of Lemma \ref{Grenzwert}]
Introduce the shorthand notation
\begin{equation*}
Z(\underline{a})\equiv Z(\sk; A, B,\underline{a})\quad \text{and}\quad
S(\underline{a})\equiv S(\sk; A, B,\underline{a}).
\end{equation*}
{}From Theorem \ref{3.2inKS1} and Lemma \ref{lem:det:neq:0} it follows that
\begin{equation*}
S(\underline{a}) = - \begin{pmatrix}\1 & 0 & 0 \end{pmatrix}
Z(\underline{a})^{-1} P_{\Ker Z(\underline{a})^\dagger}^\perp (A-\i\sk B)
\begin{pmatrix}\1 \\ 0 \\ 0 \end{pmatrix}
\end{equation*}
and
\begin{equation*}
\begin{split}
S(\underline{a}_j) & = - \begin{pmatrix}\1 & 0 & 0 \end{pmatrix}
Z(\underline{a}_j)^{-1}  (A-\i\sk B)
\begin{pmatrix}\1 \\ 0 \\ 0 \end{pmatrix}\\
& = - \begin{pmatrix}\1 & 0 & 0 \end{pmatrix} Z(\underline{a}_j)^{-1}
P_{\Ker Z(\underline{a})^\dagger}^\perp (A-\i\sk B)
\begin{pmatrix}\1 \\ 0 \\ 0 \end{pmatrix}.
\end{split}
\end{equation*}
Thus, to prove the claim it suffices to show that
\begin{equation}\label{Ziel:Grenzwert}
\lim_{j\rightarrow\infty}\begin{pmatrix}\1 & 0 & 0 \end{pmatrix}
Z(\underline{a}_j)^{-1} P_{\Ker Z(\underline{a})^\dagger}^\perp =
\begin{pmatrix}\1 & 0 & 0 \end{pmatrix} Z(\underline{a})^{-1} P_{\Ker
Z(\underline{a})^\dagger}^\perp.
\end{equation}

{}From Theorem 3.1 in \cite{KS1} it follows that all elements $z$ of $\Ker
Z(\underline{a})$ satisfy $\cP_{\cE}z=0$, where $\cP_{\cE}$ is the
orthogonal projection in $\cK$ onto $\cK_{\cE}$. Thus,
\begin{equation}\label{6.7:neu}
\begin{pmatrix}\1 & 0 & 0 \end{pmatrix} Z(\underline{a}_j)^{-1}P_{\Ker
Z(\underline{a})^\dagger}^\perp =
\begin{pmatrix}\1 & 0 & 0 \end{pmatrix} P_{\Ker Z(\underline{a})}^\perp
Z(\underline{a}_j)^{-1}P_{\Ker Z(\underline{a})^\dagger}^\perp
\end{equation}
for any $j\in\N$ and
\begin{equation}\label{6.7:neu:neu}
\begin{pmatrix}\1 & 0 & 0 \end{pmatrix} Z(\underline{a})^{-1}P_{\Ker
Z(\underline{a})^\dagger}^\perp =
\begin{pmatrix}\1 & 0 & 0 \end{pmatrix} P_{\Ker Z(\underline{a})}^\perp
Z(\underline{a})^{-1} P_{\Ker Z(\underline{a})^\dagger}^\perp.
\end{equation}
By Lemma \ref{Grenzwert:Hilfe} we have
\begin{equation*}
\lim_{j\rightarrow\infty}P_{\Ker Z(\underline{a})}^\perp
Z(\underline{a}_j)^{-1}P_{\Ker Z(\underline{a})^\dagger}^\perp = P_{\Ker
Z(\underline{a})}^\perp Z(\underline{a})^{-1}P_{\Ker
Z(\underline{a})^\dagger}^\perp.
\end{equation*}
Combining this with \eqref{6.7:neu} and \eqref{6.7:neu:neu} we obtain
\eqref{Ziel:Grenzwert}.
\end{proof}

For fixed $\sk>0$, consider
\begin{equation}\label{F:def}
F(\underline{t}) := S(\sk; A,B,\underline{a})\quad \text{with}\quad
\underline{t}=\e^{\i \underline{a} \sk}.
\end{equation}
Recall that by \eqref{11} -- \eqref{zet} the scattering matrix $S(\sk;
A,B,\underline{a})$ depends on $\underline{a}$ only through
$\underline{t}=\e^{\i \underline{a} \sk}$. The map $\underline{a}\mapsto
\e^{\i\underline{a}\sk}$ maps the set
\begin{equation*}
 \left\{\underline{a}\in\C^{|\cI|}\big|\:
\underline{a}=\{a_i\}_{i\in\cI}\; \text{with}\; 0<\Re a_i\leq
2\pi/\sk\;\text{and}\; \Im a_i > 0\; \text{for all}\; i\in\cI \right\}.
\end{equation*}
bijectively onto the polydisc $\D^{|\cI|} = \{\zeta\in\C|\: |\zeta| <
1\}^{|\cI|}$. The interval $(0,2\pi/\sk]$ is mapped onto the torus
$\T^{|\cI|} = \{\zeta\in\C|\: |\zeta| = 1\}^{|\cI|}$.

\begin{lemma}\label{lem:Hardy}
The function $F$ belongs to the Hardy class $H^p(\T^{|\cI|})$ for all
$p\in(0,\infty]$ and is inner.
\end{lemma}

\begin{remark}
We recall that an operator-valued function on a polydisc $\D^d$ is said to
be inner if it is holomorphic in $\D^d$ and takes
unitary values for almost all points of $\T^d\subset\partial(\D^d)$ (the so
called distinguished boundary of $\D^d$ \cite{Hoermander}). For $d=1$
matrix-valued inner functions are studied, e.g., in \cite{Potapov}. In
particular, an analog of the canonical factorization theorem for
matrix-valued inner functions has been proven there.
\end{remark}

\begin{proof}
{}From Proposition \ref{thm:6.1:harmo} it follows that $F$ is holomorphic
in the punctured open polydisc $\D^{|\cI|}\setminus\{0\}$. By
\eqref{S:leq:1} we have $\|F(\underline{t})\|\leq 1$ for all
$\underline{t}\in\D^{|\cI|}\setminus\{0\}$. Therefore, the Laurent
expansion of $F$ contains no terms with negative powers. Thus, $F$ is
holomorphic in $\D^{|\cI|}$.

The bound \eqref{S:leq:1} also implies that
\begin{equation*}
\sup_{r\in[0,1)}\int_{\T^{|\cI|}} \|F(r\underline{t})\|^p
d\mu(\underline{t})\leq \mu(\T^{|\cI|})
\end{equation*}
for any $p\in(0,\infty)$, where $\mu$ stands for the Haar measure on the
torus $\T^{|\cI|}$ and
\begin{equation*}
\sup_{r\in[0,1)} \sup_{\underline{t}\in\T^{|\cI|}} \|F(r\underline{t})\|
\leq 1.
\end{equation*}
For every $\underline{t}\in\T^{|\cI|}$ the operator $F(\underline{t})$ is
unitary, which means that $F$ is an inner function.
\end{proof}

\begin{proof}[Proof of Theorem \ref{thm:main:harmony}]
Since by Proposition \ref{thm:6.1:harmo} $F(\underline{t})$ is a rational
$\cB(\cK_\cE)$-valued function, it can be analytically continued as a
meromorphic function on all of $\underline{t}\in\C^{|\cI|}$. Moreover, it
is holomorphic in the polydisc $\D_{1+\epsilon}^{|\cI|}=\{\zeta\in\C|\:
|\zeta|<1+\epsilon\}$ for some $\epsilon>0$. To show this, by Hartogs'
theorem it suffices to consider the analytic continuation with respect to a
single variable $t_i\in\C$ keeping all other variables fixed. By the bound
\eqref{S:leq:1} all possible poles of this continuation lie outside a disc
$\{t_i\in\C|\;|t_i|<r\}$ with $r>1$.

In turn, this implies (see, e.g., Theorem 2.4.5 in \cite{Hoermander}) that
the Taylor series of the function $F(\underline{t})$ converges absolutely
and uniformly for all $\underline{t}\in\T^{|\cI|}$. Combining this with
Lemma \ref{Grenzwert} proves the absolute and uniform convergence of the
Fourier expansion \eqref{Fourier:exp}.

By Lemma \ref{lem:Hardy} the Fourier coefficients \eqref{fourier:coef}
satisfy $\widehat{S}_{\underline{n}}(\sk;A,B)=0$ for any $\underline{n}\in
\Z^{|\cI|}$ with $n_i<0$ for at least one $i\in\cI$.
\end{proof}

\begin{definition}\label{def:2.3}
Given a non-compact graph $\cG=(V,\cI,\cE,\partial)$ to any vertex $v\in
V=V(\cG)$ we associate the single-vertex graph
$\cG_v=(\{v\},\cI_v,\cE_v,\partial_v)$ with the following properties
\begin{itemize}
\item[(i)]{$\cI_v=\emptyset$,}
\item[(ii)]{$\partial_v(e)=v$ for all $e\in\cE_v$,}
\item[(iii)]{$|\cE_v|=\deg_{\cG}(v)$, the degree of the vertex $v$ in the graph $\cG$,}
\item[(iv)]{there is an injective map $\Psi_v:\; \cE_v\rightarrow\cE\cup\cI$ such that
$v\in\partial\circ\Psi_v(e)$ for all $e\in\cE_v$.}
\end{itemize}
\end{definition}

Since the boundary conditions are assumed to be local (see Definition
\ref{def:local}), we can consider the Laplace operator $\Delta(A_v,B_v)$ on
$L^2(\cG_v)$ associated with the boundary conditions $(A_v,B_v)$ induced by
$(A,B)$, see \eqref{permut}. By \eqref{svertex} the scattering matrix for
$\Delta(A_v, B_v )$ is given by
\begin{equation*}
S_v(\sk) = - (A_v+\i\sk B_v)^{-1} (A_v-\i\sk B_v).
\end{equation*}

Now to each walk $\bw=\{e^\prime,i_1,\ldots,i_N,e\}$ from $e^\prime\in\cE$
to $e\in\cE$ on the graph $\cG$ similar to \eqref{weight:2} we associate a
weight $W(\bw;\sk)$ by
\begin{equation}\label{weight}
W(\bw;\sk)=\e^{\i\sk\langle\underline{n}(\bw),
  \underline{a}\rangle}\,\widetilde{W}(\bw;\sk)
\end{equation}
with
\begin{equation}\label{weight:3}
\widetilde{W}(\bw;\sk)=\prod_{k=0}^{|\bw|_{\mathrm{comb}}}S_{v_{k}}(\sk)_{e_{k}^{(+)}
e_{k}^{(-)}}.
\end{equation}
Here $e_k^{(\pm)}\in\cE_{v_k}$ are defined as
\begin{equation*}
e_k^{(-)}=\begin{cases}
\Psi_{v_k}^{-1}(i_k), & \text{if}\quad 1\leq k\leq |\bw|_{\mathrm{comb}},\\
\Psi_{v_k}^{-1}(e), & \text{if}\quad k=0,
\end{cases}
\end{equation*}
and
\begin{equation*}
e_k^{(+)}=\begin{cases}
\Psi_{v_k}^{-1}(i_{k+1}), & \text{if}\quad 0\leq k\leq |\bw|_{\mathrm{comb}}-1,\\
\Psi_{v_k}^{-1}(e^\prime), & \text{if}\quad k=|\bw|_{\mathrm{comb}}+1,
\end{cases}
\end{equation*}
where the map $\Psi_v$ is defined by Definition \ref{def:2.3}. Note that
$\widetilde{W}(\bw;\sk)$ is independent of the metric properties of the
graph. Obviously, for a trivial walk $\bw=\{e^\prime, e\}$ we have
$\widetilde{W}(\bw;\sk)= S_v(\sk)_{\Psi_{v}^{-1}(e),
\Psi_{v}^{-1}(e^\prime)}$, where $v=\partial(e)=\partial(e^\prime)$.

\begin{theorem}\label{verloren}
The matrix elements of the $\underline{n}$-th Fourier coefficients
\eqref{fourier:coef} are given by the sum over the walks with score
$\underline{n}$,
\begin{equation}\label{verloren:2}
[\widehat{S}_{\underline{n}}(\sk;A,B)]_{e,e^\prime} =
\sum_{\bw\in\cW_{e,e^\prime}(\underline{n})} \widetilde{W}(\bw;\sk)
\end{equation}
if $\cW_{e,e^{\prime}}(\underline{n})$ is nonempty and
$[\widehat{S}_{\underline{n}}(\sk;A,B)]_{e,e^\prime}=0$ whenever
$\cW_{e,e^{\prime}}(\underline{n})=\emptyset$.
\end{theorem}

\begin{proof}
Obviously, it suffices to show that the $\underline{n}$-th coefficient of
the multi-dimensional Taylor expansion of the scattering matrix $S(\sk;
A,B,\underline{a})$ with respect to
$\underline{t}=\{t_i\}_{i\in\cI}\in\D^{|\cI|}$ with $t_i:=\e^{\i\sk a_i}$
coincides with the r.h.s.\ of \eqref{verloren:2}. Recall that by Theorem
\ref{3.2inKS1} and Lemma \ref{lem:det:neq:0} for all $\underline{a}\in\cA$
the scattering matrix is given by
\begin{equation}\label{S:verloren}
S(\sk; A,B,\underline{a})=-\begin{pmatrix} \1 & 0 & 0 \end{pmatrix} \left(A
X(\sk;\underline{a})+\i\sk B Y(\sk;\underline{a})\right)^{-1}(A-\i\sk B)
\begin{pmatrix} \1 \\ 0 \\ 0 \end{pmatrix},
\end{equation}
where $X(\sk;\underline{a})$ and $Y(\sk;\underline{a})$ were defined in
\eqref{zet}. Obviously,
\begin{equation}\label{verloren:3}
\begin{split}
A X(\sk;\underline{a})&+\i\sk B Y(\sk;\underline{a}) = (A+\i\sk
B)U(\sk;\underline{a}) + (A-\i\sk B)R(\sk;\underline{a})\\
&=(A+\i\sk B)\big[\1+(A+\i\sk B)^{-1}(A-\i\sk
B)R(\sk;\underline{a})U(\sk;\underline{a})^{-1}\big]U(\sk;\underline{a}),
\end{split}
\end{equation}
where
\begin{equation}\label{U:def}
U(\sk;\underline{a}):=\begin{pmatrix} \1 & 0 & 0 \\ 0 & \1 & 0 \\ 0 & 0 &
\e^{-\i\sk\underline{a}}
\end{pmatrix}
\end{equation}
and
\begin{equation*}
R(\sk;\underline{a}) := X(\sk;\underline{a}) - U(\sk;\underline{a}) =
\begin{pmatrix} 0 & 0 & 0 \\ 0 & 0 & \1 \\ 0 & \e^{\i\sk\underline{a}} & 0
\end{pmatrix}
\end{equation*}
with respect to the orthogonal decomposition \eqref{K:def}. Equation
\eqref{verloren:3} implies that
\begin{equation*}
\begin{split}
& \left(A X(\sk;\underline{a})+\i\sk B
Y(\sk;\underline{a})\right)^{-1}\\&\quad =
U(\sk;\underline{a})^{-1}\sum_{n=0}^\infty \left[-(A+\i\sk B)^{-1} (A-\i\sk
B) G H(\sk;\underline{a})\right]^n (A+\i\sk B)^{-1}
\end{split}
\end{equation*}
with
\begin{equation}\label{G:H:def}
G = \begin{pmatrix} 0 & 0 & 0 \\ 0 & 0 & \1 \\ 0 & \1 & 0
\end{pmatrix}\qquad\text{and}\qquad H(\sk;\underline{a})
= \begin{pmatrix} \1 & 0 & 0 \\ 0 & \e^{\i\sk\underline{a}} & 0 \\ 0 & 0 &
\e^{\i\sk\underline{a}}
\end{pmatrix}
\end{equation}
such that $R(\sk;\underline{a}) U(\sk;\underline{a})^{-1} = G
H(\sk;\underline{a})$. Combining this representation with
\eqref{S:verloren} we obtain
\begin{equation}\label{Reihe}
\begin{split}
& S(\sk; A,B,\underline{a})=\\ & = \sum_{n=0}^\infty
\begin{pmatrix} \1 & 0 & 0
\end{pmatrix} \left[S(\sk; A,B) G H(\sk;\underline{a}) \right]^n S(\sk; A,B)
\begin{pmatrix} \1 \\ 0 \\ 0
\end{pmatrix},
\end{split}
\end{equation}
where $S(\sk; A,B)$ is defined by \eqref{svertex}. By the unitarity of
$S(\sk; A,B)$, the series converges absolutely for all
$\underline{a}\in\cA$.

Recall that
\begin{equation*}
S(\sk; A,B) = S(\sk; CA, CB)
\end{equation*}
for every invertible $C$. It follows directly from Definition
\ref{def:local} that
\begin{equation*}
S(\sk; A,B) = \bigoplus_{v\in V(\cG)} S(\sk; A(v),B(v)).
\end{equation*}
Plugging this equality in \eqref{Reihe} proves the claim.
\end{proof}

\begin{remark}
Theorem \ref{verloren} implies that the scattering matrix of the graph $\cG$
is determined by the scattering matrices associated with  all its single
vertex subgraphs. This result can also be obtained by applying the
factorization formula \cite{KS3}.
\end{remark}

Combining Theorems \ref{thm:main:harmony} and \ref{verloren} we immediately
obtain

\begin{corollary}\label{4:prop:0}
Let $\underline{a}\in(\R_+)^{|\cI|}$ be arbitrary. For all $\sk>0$ the
scattering matrix $S(\sk;A,B,\underline{a})$ associated with the Laplacian
$\Delta(A,B,\underline{a})$ on the graph $\cG$ has an absolutely convergent
expansion in the form
\begin{equation}\label{sseries}
S(\sk;A,B,\underline{a})_{e,e^{\prime}}=\sum_{\bw\in \cW_{e,e^{\prime}}}
W(\bw; \sk) \equiv \sum_{\bw\in \cW_{e,e^{\prime}}}\widetilde{W}(\bw; \sk)
\e^{\i\sk|\bw|}.
\end{equation}
\end{corollary}

Since $\widetilde{W}(\bw;\sk)$ is independent of the metric properties of
the graph, it is natural to view \eqref{sseries} as a \emph{combinatorial
Fourier expansion} of the scattering matrix $S(\sk;A,B,\underline{a})$. We
will show that \eqref{sseries} actually coincides with the Fourier expansion
\eqref{Fourier:exp} in Theorem \ref{thm:main:harmony}.

\section{Analytic Continuation of the Scattering Matrix}\label{sec:analyt}

Recall that the scattering matrix $S(\sk;A,B,\underline{a})$ is analytic in
$\sk$ for all $\Re\sk>0$ and $\Im\sk>0$. In this section we will show that
representation \eqref{sseries} for the scattering matrix can be extended to
the complex plane.

\begin{lemma}\label{lem:beta}
There is $\beta_0>0$ such that the series
\begin{equation}\label{ssum:prime:3}
\sum_{\underline{n}\in\cN_{e,e^\prime}}
[\widehat{S}_{\underline{n}}(\sk;A,B)]_{e,e^\prime} \e^{\i\sk
\langle\underline{n},\underline{a}\rangle}
\end{equation}
converges absolutely for all $\sk\in\C$ with $\Re\sk>0$ and $\Im\sk\geq
\beta_0$. Therefore,
\begin{equation*}
\widetilde{S}(\sk; A,B,\underline{a})_{e,e^{\prime}} =
\sum_{\underline{n}\in\cN_{e,e^\prime}}
[\widehat{S}_{\underline{n}}(\sk;A,B)]_{e,e^\prime} \e^{\i\sk
\langle\underline{n},\underline{a}\rangle}
\end{equation*}
is a holomorphic function for all such $\sk\in\C$.
\end{lemma}

\begin{proof}
{}From Proposition \ref{proposition:2.3} it follows that for all
sufficiently large $\beta>0$ there is a constant $C_\beta>0$ such that the
estimate
\begin{equation*}
|S_v(\sk)_{e_1,e_2}| \leq C_\beta
\end{equation*}
holds for all $\sk\in\C$ with $\Im\sk>\beta$, any $v\in V$, and any
$e_1,e_2\in \cE_v$ (see Definition \ref{def:2.3}). Therefore, for an
arbitrary walk $\bw\in\cW_{e,e^\prime}(\underline{n})$ we obtain
\begin{equation*}
|\widetilde{W}(\bw;\sk)|  \leq C_\beta^{|\underline{n}|+1}.
\end{equation*}
Thus, from \eqref{verloren:2} it follows that
\begin{equation*}
|[\widehat{S}_{\underline{n}}(\sk;A,B)]_{e,e^\prime}| \leq
C_\beta^{|\underline{n}|+1}
 |\cW_{e,e^\prime}(\underline{n})|.
\end{equation*}
Therefore, from Lemma \ref{4:lem:2} using the identity \eqref{identity} we
obtain the estimate
\begin{equation}\label{esti}
\sum_{\substack{\underline{n}\in\cN_{e,e^\prime}\\ |\underline{n}|=N
}}|[\widehat{S}_{\underline{n}}(\sk;A,B)]_{e,e^\prime}| \leq C_\beta^{N+1}
|\cI|^{N}.
\end{equation}

Recalling the definition \eqref{amin} for $a_{\mathrm{min}}$ estimate
\eqref{esti} implies that the series
\begin{equation*}
\begin{split}
\sum_{\underline{n}\in\cN_{e,e^\prime}}
|[\widehat{S}_{\underline{n}}(\sk;A,B)]_{e,e^\prime} \e^{\i\sk
\langle\underline{n},\underline{a}\rangle}| & \leq
\sum_{\underline{n}\in\cN_{e,e^\prime}}
|[\widehat{S}_{\underline{n}}(\sk;A,B)]_{e,e^\prime}|\;
\e^{-|\underline{n}|\, \Im\sk\,
 a_{\mathrm{min}}}\\
& \leq \sum_{N=0}^\infty\e^{-N\, \Im\sk\, a_{\mathrm{min}}} \sum_{\substack{\underline{n}\in\cN_{e,e^\prime}\\
|\underline{n}|=N}} |S(\sk;\underline{n})_{e,e^{\prime}}|
\end{split}
\end{equation*}
converges for all $\sk\in\C$ with
\begin{equation*}
\Im\sk > \beta_1:=\frac{1}{a_{\mathrm{min}}} \log\left\{ C_\beta |\cI|
\right\}.
\end{equation*}
This proves the claim with $\beta_0=\max\{\beta,\beta_1\}$.
\end{proof}

The following statement is the main result of this section.

\begin{theorem}\label{ueberlap}
There is $\beta_0>0$ such that
\begin{equation}\label{ssum:prime:2}
S(\sk;A,B,\underline{a})_{e,e^{\prime}} =
\sum_{\underline{n}\in\cN_{e,e^\prime}}
[\widehat{S}_{\underline{n}}(\sk;A,B)]_{e,e^\prime} \e^{\i\sk
\langle\underline{n},\underline{a}\rangle}
\end{equation}
holds for all $\sk\in\C$ with $\Re\sk>0$ and $\Im\sk>\beta_0$.
\end{theorem}

The intuitive idea behind the proof of Theorem \ref{ueberlap} is the
observation that the series \eqref{sseries} and \eqref{ssum:prime:3} agree.
However, Theorem \ref{4:prop:0} and Lemma \ref{lem:beta} establish
convergence of these series in two disjoint sets of the complex plane.
Therefore, to prove that both series define the same analytic function we
perform a two-step analytic continuation invoking an auxiliary analytic
function of two complex variables.

\begin{proof}
Consider the $|\cE|\times|\cE|$ matrix-valued function $F(\sk_1,\sk_2)$
with matrix elements
\begin{equation*}
F(\sk_1,\sk_2)_{e,e^{\prime}}:=\sum_{\underline{n}\in\cN_{e,e^\prime}}
[\widehat{S}_{\underline{n}}(\sk_1;A,B)]_{e,e^\prime} \e^{\i\sk_2
\langle\underline{n},\underline{a}\rangle}.
\end{equation*}
By Theorem \ref{4:prop:0}
\begin{equation*}
F(\sk,\sk)_{e,e^{\prime}} = S(\sk;A,B,\underline{a})_{e,e^{\prime}}
\end{equation*}
for all $\sk=\Re \sk >0$ and by Lemma \ref{lem:beta}
\begin{equation}\label{recall}
F(\sk,\sk)_{e,e^{\prime}} = \widetilde{S}(\sk;A,B,\underline{a})_{e,e^{\prime}}
\end{equation}
for all $\sk\in\C$ with $\Re\sk > 0$ and $\Im\sk > \beta_0$, where
$\beta_0$ is defined in Lemma \ref{lem:beta}. Observe that for any
$\sk_1>0$ the function $F(\sk_1,\sk_2)$ is holomorphic in
$\sk_2\in\{\sk\in\C|\; \Re\sk>0,\: \Im\sk>0\}$. Assume that
$\Im\sk_2>\beta_0$ with $\beta_0$ defined as in Lemma \ref{lem:beta}.
Inspecting the estimates used in the proof of Lemma \ref{lem:beta} we
obtain that
\begin{equation*}
\sum_{\underline{n}\in\cN_{e,e^\prime}}
[\widehat{S}_{\underline{n}}(\sk_1;A,B)]_{e,e^\prime} \e^{\i\sk_2
\langle\underline{n},\underline{a}\rangle}
\end{equation*}
converges absolutely for all $\sk_1\in\C$ with $\Re\sk_1>0$ and
$0\leq\Im\sk_1 < \Im\sk_2+\epsilon$, where $\epsilon>0$ is sufficiently
small. Recalling \eqref{recall} completes the proof.
\end{proof}

\begin{remark}\label{rem:important}
Assume that the series
\begin{equation*}
\sum_{\underline{n}\in\cN_{e,e^\prime}}
[\widehat{S}_{\underline{n}}(\sk;A,B)]_{e,e^\prime} \e^{\i\sk
\langle\underline{n},\underline{a}\rangle}
\end{equation*}
absolutely converges in a ball $B_r(\sk_0)$ centered at $\sk_0\in\C$ with
$\Re\sk=0$ and $\Im\sk>0$. Then arguments used in the proof of Theorem
\ref{ueberlap} show that
\begin{equation*}
S(\sk;A,B,\underline{a})_{e,e^{\prime}} =
\sum_{\underline{n}\in\cN_{e,e^\prime}}
[\widehat{S}_{\underline{n}}(\sk;A,B)]_{e,e^\prime} \e^{\i\sk
\langle\underline{n},\underline{a}\rangle}
\end{equation*}
for all $\sk\in B_r(\sk_0)$.
\end{remark}

\section{The Generating Function}\label{sec:3}

In this section we prove an explicit algebraic representation for the
matrix-valued generating function $T(\beta)$ defined in equation
\eqref{statesum}. This result is formulated below as Theorem \ref{thm:3.2}.

Let $\cB$ be the canonical orthonormal basis in $\C^{|\cE|+2|\cI|}\cong
\cK=\cK_{\cE}\oplus\cK_{\cI}^{(-)}\oplus\cK_{\cI}^{(+)}$ such that any
element $h\in\cB$ is uniquely associated with some edge $j(h)\in\cI\cup\cE$.
Moreover, $j(h)\in\cE$ if $h\in\cK_{\cE}$ and $j(h)\in\cI$ if
$h\in\cK_{\cI}^{(-)}$ or $h\in\cK_{\cI}^{(+)}$. Set
\begin{equation*}
v(h) = \begin{cases} \partial(j(h)) & \text{if}\quad h\in\cK_{\cE},\\
\partial^{-}(j(h)) & \text{if}\quad h\in\cK_{\cI}^{-},\\
\partial^{+}(j(h)) & \text{if}\quad h\in\cK_{\cI}^{+}.
\end{cases}
\end{equation*}

Given a collection of matrices $\cM=\{M(v)\}_{v\in V}$ we define the
linear transformation $\mathbf{M}$ on the finite-dimensional Hilbert space $\cK$ via its
sesquilinear form
\begin{equation}\label{M:def}
\langle h_1, \mathbf{M} h_2\rangle_{\cK}=\begin{cases}
\left[M(v(h))\right]_{j(h_1),j(h_2)} & \text{if}\quad v(h_1)=v(h_2),\\
0, & \text{otherwise}.
\end{cases}
\end{equation}

For an arbitrary $\beta>0$ and every $v\in V(\cG)$ we set
\begin{equation}\label{av:bv}
A_v(\beta) := \frac{1}{2}(\1-M(v)),\qquad B_v(\beta) :=
-\frac{1}{2\beta}(\1+M(v)).
\end{equation}
Define
\begin{equation}\label{RB}
A(\beta) := \bigoplus_{v\in V} A_v(\beta),\qquad B(\beta) := \bigoplus_{v\in
V} B_v(\beta).
\end{equation}
Finally, we set
\begin{equation}\label{D}
\begin{split}
D(\beta) & = Z(\i\beta; A(\beta), B(\beta),\underline{a})\\
&= \frac{1}{2}\left(X(\i\beta;\underline{a}) +
Y(\i\beta;\underline{a})\right)-
\frac{1}{2}\mathbf{M}\left(X(\i\beta;\underline{a})-Y(\i\beta;\underline{a})\right)\\
&=\left[\1 + \mathbf{M}\, G\, H(\i\beta;\underline{a})\right]
U(\i\beta,\underline{a}).
\end{split}
\end{equation}
Here the matrix $Z(\sk;A,B,\underline{a})$ is defined in \eqref{Z:def}, the
matrices $X(\sk;\underline{a})$ and $Y(\sk;\underline{a})$ are defined in
\eqref{zet}, $U(\sk;\underline{a})$, $G$, and $H(\sk;\underline{a})$ - in
\eqref{U:def} and \eqref{G:H:def}. Writing the matrix $\mathbf{M}$ with
respect to the orthogonal decomposition \eqref{K:def} as a $3\times 3$
block-matrix
\begin{equation*}
\mathbf{M}=\begin{pmatrix} \mathbf{M}_{11} & \mathbf{M}_{12} &
\mathbf{M}_{13} \\ \mathbf{M}_{21} & \mathbf{M}_{22} & \mathbf{M}_{23} \\
\mathbf{M}_{31} & \mathbf{M}_{32} & \mathbf{M}_{33}
\end{pmatrix},
\end{equation*}
we obtain
\begin{equation*}
\1 + \mathbf{M}\, G\, H(\i\beta;\underline{a})  = \begin{pmatrix} \1 &
\mathbf{M}_{13} \e^{-\beta\underline{a}} & \mathbf{M}_{12}
\e^{-\beta\underline{a}} \\ 0 & \1+ \mathbf{M}_{23}
\e^{-\beta\underline{a}} & \mathbf{M}_{22} \e^{-\beta\underline{a}}\\
0 & \mathbf{M}_{33} \e^{-\beta\underline{a}} & \1 + \mathbf{M}_{32}
\e^{-\beta\underline{a}}
\end{pmatrix}.
\end{equation*}
Obviously, $\1 + \mathbf{M}\, G\, H(\i\beta;\underline{a})$ is an entire
matrix valued function in the complex variable $\beta$. Moreover,
\begin{equation*}
\lim_{\Re\beta\rightarrow +\infty} \det(\1+\mathbf{M}\, G\,
H(\i\beta;\underline{a})) = 1.
\end{equation*}
Thus, $\det D(\beta)$ is not identically vanishing and this in turn gives

\begin{lemma}\label{lem:6.1}
The matrix valued function $D(\beta)$ is entire in $\beta\in\C$ and its
determinant vanishes on a discrete set $\cD\subset \C$ depending on
$\underline{a}\in\R^{|\cI|}$ and the set of matrices
$\cM=\{M(v)\}_{v\in\cV}$. The set $\cD$ has no accumulation points in $\C$.
In particular, the matrix inverse $D(\beta)^{-1}$ is a meromorphic function
in $\beta\in\C$ with poles in $\cD$.
\end{lemma}

Now we turn to the main result of this article:

\begin{theorem}\label{thm:3.2}
For a given non-compact graph $\cG=(V,\cI,\cE,\partial)$ with lengths
$\underline{a}$ of the internal lines and a collection of matrices
$\cM=\{M(v)\}_{v\in V}$ at the vertices of the graph the generating function
$T(\beta)$ defined by \eqref{statesum} has an analytic extension to
$\C\setminus\cD$ and can be expressed in terms of the matrix
$D(\beta)^{-1}{\mathbf M}$ as follows
\begin{equation}\label{trep}
T(\beta)=\begin{pmatrix}\1 & 0 & 0 \end{pmatrix}D(\beta)^{-1}{\mathbf
M}\begin{pmatrix}\1 \\ 0 \\ 0 \end{pmatrix}.
\end{equation}
\end{theorem}

We turn to the proof of this theorem. First we assume that all matrices
$M(v)$ are self-adjoint. Let $A_v(\beta)$ and $B_v(\beta)$ be defined by
\eqref{av:bv}. Then$A_v(\beta) B_v(\beta)^\dagger$ is
self-adjoint, since $M(v)$ is. Observe that
\begin{equation*}
\dim \Ker\, (A_v(\beta), B_v(\beta)) = \deg(v)-\dim(\Ker\, A_v(\beta)\cap\Ker\, B_v(\beta)).
\end{equation*}
{}From \eqref{av:bv} it follows that $\Ker\, A_v(\beta)\cap\Ker\,
B_v(\beta)=\{0\}$. Therefore, the $2\deg(v)\times\deg(v)$ matrix
$(A_v(\beta), B_v(\beta))$ has maximal rank. Thus, the operator
$\Delta(A_v(\beta), B_v(\beta))$  for the single-vertex graph $\cG_v$ (see
Definition \ref{def:2.3}) is self-adjoint. The associated scattering matrix
given by \eqref{svertex} obviously satisfies the relation
\begin{equation}\label{M-S}
S(\i\beta; A_v(\beta), B_v(\beta)) = M(v).
\end{equation}

Set $B_r(\beta)=\{\sk\in\C|\; |\sk-\i\beta|<r\}$.

\begin{lemma}\label{lem:7.1}
The scattering matrix $S(\sk; A_v(\beta), B_v(\beta))$ is holomorphic for
all
\begin{equation*}
\sk\in\left( \C_+\setminus [0,\i\infty)\right)\cup B_r(\beta)
\end{equation*}
with
\begin{equation*}
r = \frac{\beta}{\|M(v)\|}.
\end{equation*}
\end{lemma}

\begin{proof}
Recalling Proposition \ref{proposition:2.3} observe that $S(\sk; A_v(\beta),
B_v(\beta))$ has a pole at $\sk=\i\varkappa$, $\varkappa\in\R_+$ if and only
if there is a $\chi\in\cL_v$ such that
\begin{equation*}
\left(\frac{1}{2}-\frac{\varkappa}{2\beta}\right) M(v) \chi =
\left(\frac{1}{2}+\frac{\varkappa}{2\beta}\right) \chi,
\end{equation*}
that is, $(\beta+\varkappa)(\beta-\varkappa)^{-1}$ is an eigenvalue of
$M(v)$. Therefore,
\begin{equation*}
\frac{\beta+\varkappa}{|\beta-\varkappa|} \leq \|M(v)\|,
\end{equation*}
which implies that the distance from the point $\i\beta$ to the closest
pole of the scattering matrix $S(\sk; A_v(\beta), B_v(\beta))$ is at least
$\beta\|M(v)\|^{-1}$.
\end{proof}

Via equations \eqref{diff:expression} and \eqref{Randbedingungen} the
matrices $A(\beta)$, $B(\beta)$ being defined by \eqref{RB} define the
self-adjoint Laplace operator $\Delta(A(\beta), B(\beta), \underline{a})$
with local boundary conditions (in the sense of Definition \ref{def:local}).

Now we choose $\beta$ so large that the series \eqref{statesum} converges.
Then, by \eqref{M-S}, the generating function $T_{e,e^\prime}(\beta)$ can
represented in the form
\begin{equation*}
T_{e,e^\prime}(\beta) = \sum_{\underline{n}\in\cN_{e,e^\prime}}
\widehat{S}_{\underline{n}}(\i\beta;A(\beta),B(\beta))_{e,e^\prime}
\e^{-\beta\langle\underline{n},\underline{a}\rangle},
\end{equation*}
where the coefficients
$\widehat{S}_{\underline{n}}(\i\beta;A(\beta),B(\beta))$ are defined by
\eqref{weight:3} and \eqref{verloren:2} with the boundary conditions
\eqref{RB}.

\begin{lemma}\label{lem:rho}
Assume that $\beta>\beta_0$ with $\beta_0$ satisfying \eqref{beta:0}. Then
there is $\rho>0$ such that the series
\begin{equation*}
\sum_{\underline{n}\in\cN_{e,e^\prime}}
[\widehat{S}_{\underline{n}}(\sk;A(\beta),B(\beta))]_{e,e^\prime}
\e^{\i\sk\langle\underline{n},\underline{a}\rangle}
\end{equation*}
converges absolutely for all $\sk\in B_\rho(\beta)$.
\end{lemma}

\begin{proof}
For an arbitrary $\epsilon>0$ choose $\rho>0$ so small that
\begin{equation*}
|S_v(\sk)_{e_1,e_2}| \leq \|M(v)\|(1+\epsilon)
\end{equation*}
for all $\sk\in B_\rho(\beta)$, all $e_1,e_2\in\cE_v$, and all $v\in V$. As
in the proof of Lemma \ref{lem:beta} for an arbitrary walk
$\bw\in\cW_{e,e^\prime}(\underline{n})$ we obtain the estimate
\begin{equation*}
|\widetilde{W}(\bw;\sk)| \leq m^{|\underline{n}|+1}
(1+\epsilon)^{|\underline{n}|+1},
\end{equation*}
where
\begin{equation*}
m:=\max_{v\in V}\|M(v)\|.
\end{equation*}
In turn, this implies the bound
\begin{equation*}
\sum_{\substack{\underline{n}\in\cN_{e,e^\prime}\\ |\underline{n}|=N}}
|[\widehat{S}_{\underline{n}}(\sk;A(\beta),B(\beta))]_{e,e^\prime}| \leq
m^{N+1} (1+\epsilon)^{N+1} |\cI|^{N}.
\end{equation*}
Therefore,
\begin{equation*}
\begin{split}
& \sum_{\underline{n}\in\cN_{e,e^\prime}}
|[\widehat{S}_{\underline{n}}(\sk;A(\beta),B(\beta))]_{e,e^\prime}
\e^{\i\sk\langle \underline{n}, \underline{a} \rangle}| \\ &
\qquad\qquad\leq \sum_{N=0}^\infty \e^{-N\; \Im\sk\;
a_{\mathrm{min}}}\sum_{\substack{\underline{n}\in\cN_{e,e^\prime}\\
 |\underline{n}|=N}} |[\widehat{S}_{\underline{n}}(\sk;A(\beta),B(\beta))]_{e,e^\prime}| \\
& \qquad\qquad\leq \sum_{N=0}^\infty \e^{-N\; \Im\sk\; a_{\mathrm{min}}}
m^{N+1} (1+\epsilon)^{N+1} |\cI|^{N}.
\end{split}
\end{equation*}
This series converges if
\begin{equation}\label{bigger}
\Im\sk > \frac{1}{a_{\mathrm{min}}}\log\left\{m (1+\epsilon) |\cI|\right\}.
\end{equation}

We claim that inequality \eqref{bigger} holds for all $\sk\in
B_\rho(\beta)$ if $\epsilon$ is chosen to be so small that
\begin{equation*}
(1+\epsilon)\e^{\beta_0-\beta} < 1,
\end{equation*}
and $\rho>0$ satisfies the inequality
\begin{equation*}
\rho < (\beta-\beta_0) - \frac{1}{a_{\mathrm{min}}} \log(1+\epsilon)
\end{equation*}
Indeed, under these assumptions for any $\sk\in B_\rho(\beta)$ we have
\begin{equation*}
\Im\sk>\beta-\rho>\beta_0+\frac{1}{a_{\mathrm{min}}} \log(1+\epsilon) >
\frac{1}{a_{\mathrm{min}}}\log\left\{m (1+\epsilon) |\cI|\right\}.
\end{equation*}
\end{proof}

\begin{proof}[Proof of Theorem \ref{thm:3.2}]
Assume the matrices $M(v)$ to be self-adjoint. Lemma \ref{lem:rho} and
Remark \ref{rem:important} imply that there is $\rho>0$ such that
\begin{equation*}
\sum_{\underline{n}\in\cN_{e,e^\prime}}
[\widehat{S}_{\underline{n}}(\sk;A(\beta),B(\beta))]_{e,e^\prime}
\e^{\i\sk\langle\underline{n},\underline{a}\rangle}=S(\sk;
A(\beta),B(\beta);\underline{a})
\end{equation*}
holds for all $\sk\in B_\rho(\beta)$. Thus, the generating function
$T(\beta)$ can be expressed in terms of the scattering matrix,
\begin{equation}\label{TS}
T(\beta)=S(\i\beta; A(\beta),B(\beta);\underline{a}).
\end{equation}
In turn, the scattering matrix can be calculated by means of Theorem
\ref{3.2inKS1}. Obviously,
\begin{equation*}
Z(\i\beta;A(\beta),B(\beta);\underline{a}) = D(\beta)
\end{equation*}
and
\begin{equation*}
A(\beta)+\beta B(\beta)=-\mathbf{M}.
\end{equation*}
Thus, \eqref{trep} follows from \eqref{S-matrix}.

Now we relax the assumption on the self-adjointness of the matrices $M(v)$.
Obviously, the r.h.s.\ of \eqref{trep} is a rational function with respect
to the entries of the matrix $\mathbf{M}$. Since $\det D(\beta)$ does not
vanish identically, we obtain the claim.
\end{proof}

\begin{remark}
Relation \eqref{TS} combined with the factorization
formula for the scattering matrix on the graph \cite{KS3} allows to
determine the generating function $T(\beta)$ of walks on the graph $\cG$ in
terms of the the generating functions associated with subgraphs of $\cG$.
\end{remark}

\begin{remark}\label{directproof}
There is a direct way to establish \eqref{trep}. Indeed, observe that by
\eqref{D} one has
\begin{equation}\label{trepnew}
\begin{pmatrix}\1 & 0 & 0 \end{pmatrix}D(\beta)^{-1}{\mathbf
M}\begin{pmatrix}\1 \\ 0 \\ 0 \end{pmatrix}=
\begin{pmatrix}\1 & 0 & 0 \end{pmatrix}
\left[\1 + \mathbf{M}\, G\, H(\i\beta;\underline{a})\right]^{-1} \mathbf{M}
\begin{pmatrix}\1 \\ 0 \\ 0 \end{pmatrix}.
\end{equation}
The matrix $G$ performs the ``jump'' from one boundary vertex of an
internal line to the other. A simple calculation shows that if
$\Re\beta>0$, then
\begin{equation*}
\|GH(\i\beta;\underline{a})\| = \e^{-\Re\beta\;
a_{\mathrm{min}}},\quad\text{where}\quad a_{\mathrm{min}}=\min_{i\in\cI}
a_i.
\end{equation*}
Therefore, the series expansion of $\left[\1 + \mathbf{M}\, G\,
H(\i\beta;\underline{a})\right]^{-1}$ converges absolutely for all
$\beta\in\C$ with sufficiently large $\Re\beta>0$. The expression
\eqref{trepnew} coincides with the series \eqref{statesum}.
\end{remark}

This observation gives rise to the following generalization, where the
penalty vector depends on the direction in which a given edge is traversed
by a walk. Let $\underline{a}=\{a_i\}_{\i\in\cI}$ and
$\underline{b}=\{b_i\}_{\i\in\cI}$ be two arbitrary penalty vectors. Set
\begin{equation*}
\widehat{H}(\sk;\underline{a},\underline{b}) = \begin{pmatrix} \1 & 0 & 0 \\
0 & \e^{\i\sk\underline{a}} & 0 \\ 0 & 0 & \e^{\i\sk\underline{b}}
\end{pmatrix}
\end{equation*}
such that
$\widehat{H}(\sk;\underline{a},\underline{a})=H(\sk;\underline{a})$. Define
now
\begin{equation}\label{trepsuper}
\widehat{T}(\beta)=\begin{pmatrix}\1 & 0 & 0 \end{pmatrix} \left[\1 +
\mathbf{M}\, G\,
\widehat{H}(\i\beta;\underline{a},\underline{b})\right]^{-1}\mathbf{M}
\begin{pmatrix}\1 \\ 0 \\ 0 \end{pmatrix}.
\end{equation}
For any nontrivial walk $\bw=\{e^\prime,i_1,\ldots,i_N,e\}$ we set
\begin{equation*}
c_{i_k}=\begin{cases}a_{i_k}&\text{if the walk traverses the
    edge}\; i_k\in\cI\;\text{in the direction}\\ &
    \qquad\qquad\text{from the terminal to the initial vertex,}\\
b_{i_k}&\text{if the walk traverses the
    edge}\; i_k\in\cI\;\text{in the direction}\\ &
    \qquad\qquad\text{from the initial to the terminal vertex,}
\end{cases}
\end{equation*}
where $k\in\{1,\ldots,|\bw|_{\mathrm{comb}}\}$.

Using the arguments presented above one can easily prove the following
statement.

\begin{theorem}\label{genrep}
For all $\beta\in\C$ with $\Re \beta$ being sufficiently large the function
$\widehat{T}(\beta)$ equals the generating function defined by the series
\eqref{statesum} with $a_{i_k}$ being replaced by $c_{i_k}$.
\end{theorem}

\section{Random Walks on Graphs}\label{sec:random}

In this section we define random walks on a non-compact graph $\cG$ endowed
with the metric structure given by a penalty vector $\underline{a}$. Assume
that the matrices $M(v)$ are stochastic, that is, all their entries are
nonnegative and satisfy
\begin{equation*}
\sum_{k_1} \left[M(v)\right]_{k_1, k_2} = 1\quad\text{for any edge}\:
k_2\in\cI\cup\cE\:\text{incident with the vertex}\; v,
\end{equation*}
where the sum is taken over all edges $k_1\in\cI\cup\cE$ incident with the
vertex $v$. The external lines of the graph will be interpreted as initial
or final states of the walk, the internal lines as intermediate states.

Take an arbitrary external line $e\in\cE$ and consider a sequence
$\{X\}_{n=0}^N$ of random variables with values in the set $\cI\cup\cE$
determined by the following rule. Set $X_0=e$. Let $v_0=\partial(e)$.
Choose randomly an element $j_1$ of $S(v_0)$ with probability
$M(v_0)_{j_1,e}$. Set $X_1=j_1$. If $j_1\in\cE$, then $N=2$ and the
sequence is completed. If $j_1\in\cI$, then take $v_1\in\partial(j_1)$,
$v_1\neq v_0$. Choose randomly an element $j_2$ of $S(v_1)$ with probability
$M(v_1)_{j_2,j_1}$ and set $X_2=j_2$. If $j_2\in\cE$, then $N=3$ and the
sequence is completed. Otherwise proceed inductively. Finally, we obtain
finite of infinite sequence of random variables. If $N<\infty$, then
$\{X\}_{n=0}^N$ is a walk in the sense of Section \ref{sec:2}. We call this
sequence a \emph{random walk} on the graph $\cG$ from $e\in\cE$ to
$e^\prime=X_N\in\cE$.

The generating function of random walks from $e\in\cE$ to $e^\prime\in\cE$
is defined by equation \eqref{statesum}. Obviously, it is monotone with
respect to $\beta$,
\begin{equation*}
T_{e,e^{\prime}}(\beta)\leq T_{e,e^{\prime}}(\beta^{\prime})
\end{equation*}
for $\beta\geq\beta^{\prime}$. If $\cW_{e,e^\prime}$ contains
at least one nontrivial walk, then
$T_{e,e^{\prime}}(\beta)$ is strictly monotone with respect to $\beta$,
\begin{equation*}
T_{e,e^{\prime}}(\beta)< T_{e,e^{\prime}}(\beta^{\prime})
\end{equation*}
for $\beta>\beta^{\prime}$.

Recall that the stochastic matrix $M(v)$ is said to be regular if it is
ergodic, i.e., if there is a natural number $k$ such that the $k$-th power
$M(v)^k$ of the matrix $M(v)$ has strictly positive matrix entries.

\begin{lemma}
Let $\cG$ be a non-compact connected graph. Assume that each $M(v)$ is a
regular stochastic matrix. If in addition all diagonal elements of each
$M(v)$ are strictly positive, then all matrix elements of $T(\beta)$ are
strictly positive for all sufficiently large $\beta>0$.
\end{lemma}

\begin{proof}
Connectedness of $\cG$ implies that all $\cW_{e,e^{\prime}}$ are non-empty.
Given $e$ and $e^{\prime}$ for $T_{e,e^{\prime}}(\beta) >0$ to hold it is
necessary and sufficient that there is at least one walk
$\bw\in\cW_{e,e^{\prime}}$ with $W(\bw) > 0$. For the last condition to hold
it is in turn sufficient that all matrices $M(v)$ are ergodic and their
diagonal elements are strictly positive.
\end{proof}

In the remainder of this section we will discuss several examples and
introduce some mean values associated with random walks on the graph $\cG$.
These mean values are related to the generating function and its derivative
evaluated at $\beta=0$. However, for $\beta>0$ the generating function
$T_{e,e^{\prime}}(\beta)$ can be interpreted as a partition function (see,
e.g., \cite{Ruelle}) with $\beta$ being the inverse temperature. The role of
the statistical ensemble is played here by the set $\cW_{e,e^\prime}(\cM)$
of all relevant walks from $e^\prime$ to $e$.

\subsection*{1.}

We leave it to the reader to verify that the mean length of a random walk from
$e^\prime\in\cE$ to $e\in\cE$ is given by
\begin{equation}\label{length}
\langle |\bw|\rangle = - \left.\frac{d}{d\beta} \log T_{e,
e^{\prime}}(\beta)\right|_{\beta=0}= - \left. T_{e,
e^{\prime}}(\beta)^{-1}\frac{d}{d\beta} T_{e,
e^{\prime}}(\beta)\right|_{\beta=0}.
\end{equation}
The r.h.s.\ of \eqref{length} can be calculated by means of Theorem
\ref{thm:3.2}. In the thermodynamic setting (i.e., for $\beta>0$) the
quantity
\begin{equation*}
- \frac{d}{d\beta} \log T_{e, e^{\prime}}(\beta)= -  T_{e,
e^{\prime}}(\beta)^{-1}\frac{d}{d\beta} T_{e, e^{\prime}}(\beta)
\end{equation*}
corresponds to the ``mean length'' at the temperature $\beta^{-1}$. In the
following examples we will consider probabilistic ($\beta=0$) and
thermodynamic ($\beta>0$) means on equal ground.

\subsection*{2.}

As another example we consider the following situation. We say that a
walker entering a vertex $v$ from the edge $k$ and leaving through the edge
$j$ experiences a transition from $k$ to $j$ at $v$. Now fix a vertex
$v_{0}\in V$ and edges $j_{0}, k_{0}\in\cG_{v_{0}}$ satisfying the
inequality $M(v_{0})_{j_{0},k_{0}}
> 0$. We set
\begin{equation}\label{ml}
M(v_{0};\lambda)_{jk}=\left\{\begin{array}{cc}
e^{-\lambda}M(v_{0})_{jk}\quad&\mbox{if}\quad j=j_{0},\,k=k_{0}\\
M(v_{0})_{jk} \quad&\mbox{otherwise}
\end{array}
\right.
\end{equation}
with an arbitrary $\lambda >0$. Note that $M(v_{0};\lambda)$ is ergodic if
$M(v_{0})$ is, but of course not stochastic. Now replacing $M(v_{0})$ by
$M(v_{0};\lambda)$ while leaving all other $M(v)$ in the collection
$\{M(v)\}_{v\in V}$ unchanged, consider the matrix $\textbf{M}(\lambda)$
defined by \eqref{M:def}. Obviously, $\mathbf{M}(0)=\mathbf{M}$.  Further,
similar to \eqref{D}, we introduce the matrix $D(\beta;\lambda)$
\begin{equation*}
D(\beta;\lambda) = \frac{1}{2}\left(X(\i\beta;\underline{a}) +
Y(\i\beta;\underline{a})\right)-
\frac{1}{2}\mathbf{M}(\lambda)\left(X(\i\beta;\underline{a})-Y(\i\beta;\underline{a})\right),
\end{equation*}
and define the generating function $T(\beta;\lambda)$ in analogy with
\eqref{statesum} by
\begin{equation*}
T_{e,e^{\prime}}(\beta;\lambda) = \sum_{\bw\in\cW_{e,e^{\prime}}}
W(\bw;\lambda) \e^{-\beta|\bw|}
\end{equation*}
with
\begin{equation*}
W(\bw;\lambda)=\prod_{k=0}^{|\bw|_{\mathrm{comb}}}\left[M({v_{k}};\lambda)\right]_{e_{k}^{(+)}
e_{k}^{(-)}}.
\end{equation*}
Obviously, Theorem \ref{thm:3.2} remains valid for $T(\beta;\lambda)$ such
that
\begin{equation}\label{lmu}
T(\beta;\lambda)=\begin{pmatrix}\1 & 0 & 0
\end{pmatrix}D(\beta;\lambda)^{-1}\mathbf{M}(\lambda)\begin{pmatrix}\1 \\ 0 \\ 0
\end{pmatrix}.
\end{equation}
Observe that if $D(\beta)$ is invertible for a given $\beta$ then
$D(\beta;\lambda)$ is also invertible for the same $\beta$ and all
sufficiently small $\lambda > 0$. We, obviously, have
\begin{equation*}
-\frac{d}{d\lambda}W(\bw;\lambda)\Big|_{\lambda=0}=
n_{v_{0},j_{0},k_{0}}(\bw)
 W(\bw; 0)=n_{v_{0},j_{0},k_{0}}(\bw)
 W(\bw),
\end{equation*}
where $n_{v_{0},j_{0},k_{0}}(\bw)\geq 0$ is the number of times a walker
experiences a transition from $k_{0}$ to $j_{0}$ at the vertex $v_{0}$
along a given walk $\bw\in\cW_{e,e^{\prime}}(\cM)$.

Consider the quantity
\begin{equation}\label{mean2}
\begin{aligned}
\langle n_{v_{0},j_{0},k_{0}}^{e,e^\prime}\rangle(\beta)
=&-\frac{d}{d\lambda}\log T_{e,e^{\prime}}(\beta;\lambda)\Big|_{\lambda=0}\\
=&-T_{e,e^{\prime}}(\beta)^{-1}\frac{d}{d\lambda}
T_{e,e^{\prime}}(\beta;\lambda)\Big|_{\lambda=0}.
\end{aligned}
\end{equation}
It is easy to verify that $\langle
n_{v_{0},j_{0},k_{0}}^{e,e^\prime}\rangle(\beta)$ is the mean number of
times a random walk from $e^{\prime}\in\cE$ to $e\in\cE$ experiences a
transition from $k_{0}$ to $j_{0}$ at the vertex $v_{0}$.

Using Theorem \ref{thm:3.2} the derivative in \eqref{mean2} can be
calculated in a rather simple way:
\begin{equation*}
\begin{aligned}
\frac{d}{d\lambda}T(\beta;\lambda)\Big|_{\lambda=0} = & -\begin{pmatrix}\1
& 0 & 0
\end{pmatrix}
\frac{d}{d\lambda}
\Big(D(\beta;\lambda)^{-1} \mathbf{M}(\lambda)\Big)\Big|_{\lambda=0}\begin{pmatrix}\1 \\ 0 \\
0
\end{pmatrix}\\
=&\begin{pmatrix}\1 & 0 & 0
\end{pmatrix} \Big(D(\beta)^{-1}
\frac{d}{d\lambda} D(\beta;\lambda)\Big|_{\lambda=0}
D(\beta)^{-1}\mathbf {M}\Big)\begin{pmatrix}\1 \\ 0 \\
0
\end{pmatrix}\\
&\quad-\begin{pmatrix}\1 & 0 & 0
\end{pmatrix} D(\beta)^{-1}
\frac{d}{d\lambda}\mathbf{M}(\lambda)\Big|_{\lambda=0}\begin{pmatrix}\1 \\ 0 \\
0
\end{pmatrix}.
\end{aligned}
\end{equation*}
Now set
\begin{equation*}
-\frac{d}{d\lambda}\mathbf{M}(\lambda)\Big|_{\lambda=0}
=\mathbf{M}(v_{0},j_{0},k_{0})
\end{equation*}
such that
\begin{equation*}
\frac{d}{d\lambda}D(\beta;\lambda)\Big|_{\lambda=0}= \frac{1}{2}{\mathbf
M}(v_{0},j_{0},k_{0})\Big(X(\i\beta;\underline{a})-Y(\i\beta;\underline{a})\Big),
\end{equation*}
where the matrices $X$ and $Y$ are defined in \eqref{zet}. Therefore,
\begin{equation*}
\begin{aligned}
\frac{d}{d\lambda}T(\beta;\lambda)\Big|_{\lambda=0} = &\;\frac{1}{2}
\begin{pmatrix}\1 & 0 & 0
\end{pmatrix} \Big(D(\beta)^{-1}\big(X(\i\beta;\underline{a}) -
Y(\i\beta;\underline{a})\big)\\ &\qquad\cdot \mathbf {M}(v_{0},j_{0},k_{0})
D(\beta)^{-1} \mathbf{M}\Big)\begin{pmatrix}\1 \\ 0 \\
0
\end{pmatrix}\\
& + \begin{pmatrix}\1 & 0 & 0
\end{pmatrix} D(\beta)^{-1}\mathbf{M}(v_{0},j_{0},k_{0})\begin{pmatrix}\1 \\ 0 \\
0
\end{pmatrix}.
\end{aligned}
\end{equation*}
Thus, only the knowledge of the inverse $D(\beta)^{-1}$ is necessary to
determine $\frac{d}{d\lambda}T(\beta;\lambda)\Big|_{\lambda=0}$. Note that
only one matrix element of $\mathbf{M}(v_{0},j_{0},k_{0})$ is non-vanishing.

The quantities
\begin{equation*}
\begin{aligned}
\langle n_{v_{0},j_{0},\bullet}^{e,e^\prime}\rangle(\beta)&=\sum_{k_{0}}
\langle n_{v_{0},j_{0},k_{0}}^{e,e^\prime}\rangle(\beta)\\
\langle n_{v_{0},\bullet,k_{0}}\rangle(\beta)&=\sum_{j_{0}} \langle
n_{v_{0},j_{0},k_{0}}\rangle(\beta)
\end{aligned}
\end{equation*}
are related to the mean values for the probability that the vertex $v_{0}$
is entered -- during a walk from $e^{\prime}$ to $e$ -- via $j_{0}\in
\cS(v_{0})$ or left via $k_{0}\in \cS(v_{0})$, respectively. Therefore,
\begin{equation*}
\begin{split}
\langle n_{v_{0}}^{e,e^\prime}\rangle(\beta) &
=\sum_{j_{0},k_{0}\in\cS(v_{0})} \langle
n_{v_{0},j_{0},k_{0}}^{e,e^\prime}\rangle(\beta)\\ &
=\sum_{j_{0}\in\cS(v_{0})} \langle
n_{v_{0},j_{0},\bullet}^{e,e^\prime}\rangle(\beta)\\
& =\sum_{k_{0}\in\cS(v_{0})} \langle
n_{v_{0},\bullet,k_{0}}^{e,e^\prime}\rangle(\beta)
\end{split}
\end{equation*}
is the mean number of times the vertex $v_{0}$ is
visited during random walks from $e^{\prime}$ to $e$. Similarly,
\begin{equation}\label{mean32}
\langle \bar{n}_{v_{0}}^{e,e^\prime}\rangle(\beta)
=\sum_{j_{0}\in\cS(v_{0})} \langle
n_{v_{0},j_{0},j_{0}}^{e,e^\prime}\rangle(\beta)
\end{equation}
is the mean number of times the vertex $v_{0}$ is
entered and left through the same edge during a walk from $e^{\prime}$ to
$e$.

Assume now that for given $e^{\prime}\in\cE$ we have
$T_{e,e^{\prime}}(\beta) > 0$ for all $e\in\cE$. Set
\begin{equation*}
T_{\bullet e^{\prime}}(\beta)=\sum_{e\in\cE} T_{e,e^{\prime}}(\beta).
\end{equation*}
Then, the value of the quantity
\begin{equation}
\label{mean33} \langle n_{v_{0}}^{\bullet e^{\prime}}\rangle(\beta)
=\sum_{e\in\cE}\langle n_{v_{0}}^{e e^{\prime}}\rangle(\beta)
\frac{T_{e,e^{\prime}}(\beta)} {T_{\bullet e^{\prime}}(\beta)}
\end{equation}
gives the mean number of visits at the vertex $v_{0}$ for
random walks starting at $e^{\prime}\in\cE$. Similarly, if for given
$e\in\cE$ we have $T_{e,e^{\prime}}(\beta) > 0$ for all $e^{\prime}\in\cE$
we set
\begin{equation*}
T_{e \bullet}(\beta)=\sum_{e^{\prime}\in\cE} T_{e,e^{\prime}}(\beta).
\end{equation*}
The quantity
\begin{equation}\label{mean34}
\langle n_{v_{0}}^{e\bullet}\rangle(\beta) =\sum_{e^{\prime}\in\cE}\langle
n_{v_{0}}^{e e^{\prime}}\rangle(\beta) \frac{T_{e,e^{\prime}}(\beta)}
{T_{e\bullet}(\beta)}
\end{equation}
is the mean number of visits of the vertex $v_{0}$ for walks
ending at $e\in\cE$. With
\begin{equation*}
T_{\bullet\bullet}(\beta)=\sum_{e,e^{\prime}\in\cE} T_{e,e^{\prime}}(\beta)=
\sum_{e^{\prime}\in\cE} T_{\bullet e^{\prime}}(\beta)= \sum_{e\in\cE}
T_{e\bullet}(\beta)
\end{equation*}
consider the quantity
\begin{equation*}
\begin{aligned}
\langle n_{v_{0}}^{\bullet \bullet}\rangle(\beta)& = \sum_{e\in\cE}\langle
n_{v_{0}}^{e \bullet}\rangle(\beta) \frac{T_{e\bullet}(\beta)} {T_{\bullet
\bullet}(\beta)}=\sum_{e^{\prime}\in\cE}\langle n_{v_{0}}^{\bullet
e^{\prime}} \rangle(\beta)\frac{T_{\bullet e^{\prime}}(\beta)}
{T_{\bullet \bullet}(\beta)}\\
&=\sum_{e,e^{\prime}\in\cE}\langle n_{v_{0}}^{e e^{\prime}}
\rangle(\beta)\frac{T_{e,e^{\prime}}(\beta)} {T_{\bullet \bullet}(\beta)}.
\end{aligned}
\end{equation*}
Obviously, $\langle n_{v_{0}}^{\bullet \bullet}\rangle(\beta)$ is the mean
number a random walk in
$\cW(\cG)=\cup_{e,e^{\prime}}\cW_{e,e^{\prime}\in\cE}(\cG)$ visits the
vertex $v_{0}$. Therefore,
\begin{equation*}
\sum_{v_{0}\in V}\langle n_{v_{0}}^{\bullet \bullet} \rangle(\beta)\geq 1
\end{equation*}
is the mean number of vertices visited during a random walk.

\subsection*{3}

As a final example we consider the mean number $\langle
n_{i_{0}}^{ee^{\prime}}\rangle(\beta)$ any internal line $i_{0}\in\cI$ is
traversed (in either direction) by a random walk from $e^{\prime}$ to $e$.
For this replace $a_{i_{0}}$ by $a_{i_{0}} \e^\mu$ while keeping all other
$a_{i}$ fixed and set
\begin{equation*}
\underline{a}(i_0,\mu)=\{a_i(i_0,\mu)\}_{i\in\cI}\quad\text{with}\quad a_i(i_0,\mu)=
\begin{cases}a_i, & \text{if}\; i\neq i_0,\\
a_{i_{0}} \e^\mu, & \text{if}\; i = i_0.
\end{cases}
\end{equation*}
Denote by $T(\beta;\mu)$ the resulting generating function. Then
\begin{equation*}
\begin{aligned}
\langle n_{i_{0}}\rangle_{e e^{\prime}}(\beta)
=&-\frac{1}{\beta}\frac{d}{d\mu}\log T_{e,e^{\prime}}(\beta;\mu)
\Big|_{\substack{\mu=0}}\\
=&-T_{e,e^{\prime}}(\beta)^{-1}\frac{d}{d\mu}
T_{e,e^{\prime}}(\beta;\mu)\Big|_{\substack{\mu=0}}.
\end{aligned}
\end{equation*}
The derivative of the generating function with respect to $\mu$ can be
calculated by means of Theorem \ref{thm:3.2}, thus yielding,
\begin{equation*}
\begin{aligned}
\frac{d}{d\mu} T(\beta;\mu)\Big|_{\mu=0} = &\frac{1}{2} \begin{pmatrix}\1 &
0 & 0
\end{pmatrix} \Big(D(\beta)^{-1} \big((\1-\mathbf{M})\frac{d}{d\mu}
X(\i\beta;\underline{a}(i_0,\mu))\Big|_{\mu=0}\\
&+(\1+\mathbf{M})\frac{d}{d\mu}
Y(\i\beta;\underline{a}(i_0,\mu))\Big|_{\mu=0}\;\big)
D(\beta)^{-1}\mathbf{M}\Big)\begin{pmatrix}\1 \\ 0 \\ 0
\end{pmatrix}.
\end{aligned}
\end{equation*}

Similar to the discussion of the mean
number of vertices visited during a random walk we introduce the quantities
\begin{equation*}
\begin{aligned}
\langle n_{i_{0}}^{\bullet e^{\prime}}\rangle(\beta)& =
\sum_{e\in\cE}\langle n_{i_{0}}^{e e^{\prime}}\rangle(\beta)
\frac{T_{e,e^{\prime}}(\beta)}
{T_{\bullet e^{\prime}}(\beta)},\\
\langle n_{i_{0}}^{e\bullet}\rangle(\beta) & = \sum_{e^{\prime}\in\cE}
\langle n_{i_{0}}^{e e^{\prime}}\rangle(\beta)
\frac{T_{e,e^{\prime}}(\beta)}
{T_{e\bullet}(\beta)},\\
\langle n_{i_{0}}^{\bullet\bullet}\rangle(\beta) & =
\sum_{e,e^{\prime}\in\cE} \langle n_{i_{0}}^{e
e^{\prime}}\rangle(\beta)\frac{T_{e,e^{\prime}}(\beta)}
{T_{\bullet\bullet}(\beta)}.
\end{aligned}
\end{equation*}
Thus, $\langle n_{i_{0}}^{\bullet e^{\prime}}\rangle(\beta)$ is the mean
number of times the internal line $i_{0}\in\cI$ is traversed by a random
walk starting at $e^{\prime}\in\cE$, $\langle
n_{i_{0}}^{e\bullet}\rangle(\beta)$ the mean number the internal line
$i_{0}\in\cI$ is traversed by a random walk ending at $e\in\cE$. The
quantity $\langle n_{i_{0}}^{\bullet\bullet}\rangle(\beta)$ is the mean
number of times the internal line $i_{0}\in\cI$ is traversed by any random
walk.

\appendix
\section*{Appendix. Random Walks on Vertices}\label{sec:app}
\setcounter{equation}{0}
\renewcommand{\theequation}{A.\arabic{equation}}

Here we will relate the customary notion of random walks on graphs (see,
e.g., \cite{Aldous:Fill} or \cite{Woess}) to random walks considered in the
present work. Recall that the customary notion of random walks on graphs is
given by a Markov chain with vertices as states. The transition matrix $P$
indexed by the vertices has a non-vanishing entry only if the corresponding
vertices are adjacent.

Consider a graph
$\cG^\prime=\cG^\prime(V^\prime,\cI^\prime,\emptyset,\partial^\prime)$ with
no external lines. Let $P: V^{\prime}\times V^{\prime} \longrightarrow
\R_+\cup\{0\}$ be a nearest neighbor transition matrix, i.e.,
\begin{equation}\label{markov1}
\sum_{v^{\prime}\in
V^{\prime}}P(v^{\prime},v)=1\quad\mbox{for any}\quad v\in V
\end{equation}
(we read from right to left) and $P(v^{\prime},v)>0$ occurs only if $v$ and
$v^\prime$ are adjacent.

Pick an arbitrary vertex in $\cG^{\prime}$ which we denote by $v_{\infty}$.
Let $V_{v_\infty}\subset V^\prime$ be the set of all vertices adjacent to
$v_\infty$, $\cI_{v_\infty}$ the set of the internal lines $i\in\cI$
incident with $v_{\infty}$. For any $i\in\cI_{v_\infty}$ let $v_i\in
V_{v_\infty}$ be the vertex adjacent to $v_\infty$ by $i\in\cI$, that is,
\begin{equation*}
\text{either}\quad \partial^\prime(i) =
(v_\infty,v_i)\quad\text{or}\quad\partial^\prime(i) = (v_i,v_\infty)
\end{equation*}
Now replace every edge $i\in\cI_{v_\infty}$ by the external line $e$
incident with  the vertex $v_i$. Denote the set of all external lines by
$\cE$ and define the boundary operator
\begin{equation*}
\partial(j) =\begin{cases}\partial^\prime(j), & \text{if}\quad j\in\cI^\prime,\\
v, & \text{if}\quad j\in\cE. \end{cases}
\end{equation*}
Thus, we have constructed a non-compact graph $\cG(V,\cI,\cE,\partial)$ with
$V=V^\prime\setminus v_\infty$ and $\cI=\cI^\prime\setminus
\cI_{v_\infty}$. Obviously, the degree of any vertex $v\in V$ being
calculated for the graphs $\cG^\prime$ and $\cG$ is equal.

Let $\cS(v)$ be the star graph of the vertex $v\in V$, that is, the set of
all edges $j\in\cI\cup\cE$ which are incident with the vertex $v$. Given a
matrix $P$ and $v\in V$ we define the $\deg(v)\times \deg(v)$ matrix $M(v)$
with entries $M(v)_{ij}$, $i,j\in\cS(v)$ as follows:
\begin{equation*}
0\leq M(v)_{ij}=\begin{cases} P(v^\prime,v)\;\text{with}\;
v^\prime\in\partial(i),\; v^\prime\neq v, &\quad\mbox{for}
\quad i\in\cS(v)\setminus\cE,\\
P(v_{\infty},\partial(e))&\quad\mbox{for}\quad i=e\in\cS(v)\cap\cE.
\end{cases}
\end{equation*}
In particular, the matrix element $M(v)_{ij}$ is independent of $j$ and by
\eqref{markov1}
\begin{equation}\label{markov2}
\sum_{i\in\cL(v)} M(v)_{ij}=1\quad\mbox{for all}\quad
j\in\cL(v).
\end{equation}

A converse construction is also possible. Assume that a non-compact graph
$\cG$ has $\cE\neq\emptyset$ and any two vertices of the graph are adjacent
by no more than one internal line. Further, assume that all matrix entries
$0\le M(v)_{ij}$ are independent of $j$, that is, in each matrix $M(v)$ all
columns are equal, and the equality \eqref{markov2} holds. Consider the
graph $\cG^{\prime}$ without external lines obtained from $\cG$ by
replacing each external line $e$ by an internal incident with an additional
vertex $v_\infty$ such that its vertex set $V^\prime = V\cup\{v_{\infty}\}$.
Now for any $v^{\prime},v\in V^{\prime}$ we set
\begin{equation}
\label{markov3} P(v^{\prime},v)=\begin{cases}
M(v)_{ij}&\quad\mbox{for}\quad i\in\cI:\; v^\prime,v\in\partial(i),\\
M(v)_{ij}&\quad\mbox{for}\quad i\in\cE, v^{\prime}=v_{\infty}, v=\partial(e),\\
{|\cE|}^{-1} & \quad \mbox{for}\quad v=v_{\infty},\,v^\prime\in V_{v_\infty},\\
0&\quad\mbox{otherwise.}
\end{cases}
\end{equation}
Then $P$ is a nearest neighbor transfer matrix.


\end{document}